\documentclass[10pt,twoside,german]{amsart}
\advance\oddsidemargin by -1.0cm
\advance\evensidemargin by -1.0cm
\advance\topmargin by -1.0cm 
\usepackage{amssymb}
\usepackage{babel}
\usepackage{amstext}
\input{diagrams}
\usepackage{amscd}   
\usepackage{epsfig}  
\usepackage{rotating}
\let\mathg\mathfrak
\theoremstyle{plain}
\newtheorem{kor}{Korollar}
\newtheorem{lem}{Lemma}
\newtheorem{thm}{Satz}
\newtheorem{prop}{Proposition}
\theoremstyle{definition}
\newtheorem{exa}{Beispiel}
\newtheorem{NB}{Bemerkung}
\newtheorem{dfn}{Definition}
\newtheorem*{thank}{Danksagung}
\newcommand{\bdm}{\begin{displaymath}}
\newcommand{\edm}{\end{displaymath}}
\newcommand{\ba}[1]{\begin{array}{#1}}
\newcommand{\ea}{\end{array}}
\newcommand{\bea}[1][]{\begin{eqnarray#1}}
\newcommand{\eea}[1][]{\end{eqnarray#1}}
\newcommand{\btab}{\begin{tabular}}
\newcommand{\etab}{\end{tabular}}

\newcommand{\x}{\times}
\newcommand{\op}{\oplus}
\newcommand{\ox}{\otimes}
\newcommand{\bigop}{\bigoplus}

\newcommand{\ra}{\rightarrow}
\newcommand{\lra}{\longrightarrow}

\newcommand{\lmapsto}{\longmapsto}
\newcommand{\inj}{\hookrightarrow}

\newcommand{\qqs}{\forall}
\newcommand{\sla}[1]{\ensuremath{\not{\!#1}}}

\newcommand{\rk}{\ensuremath{\mathrm{rk}\,}}
\newcommand{\codim}{\ensuremath{\mathrm{codim}\,}}
\newcommand{\id}{\ensuremath{\mathrm{id}}}
\newcommand{\restr}{\arrowvert}
\newcommand{\del}{\partial}
\newcommand{\lgb}{\lbrack\!\lbrack}
\newcommand{\rgb}{\rbrack\!\rbrack}
\newcommand{\ad}{\ensuremath{\mathrm{ad}}}
\newcommand{\trdeg}{\ensuremath{\mathrm{tr\, deg}\,}}
\newcommand{\C}{\ensuremath{\mathbf{C}}}
\newcommand{\R}{\ensuremath{\mathbf{R}}}
\newcommand{\K}{\ensuremath{\mathbf{K}}}

\newcommand{\Z}{\ensuremath{\mathbf{Z}}}
\newcommand{\Q}{\ensuremath{\mathbf{Q}}}
 
\newcommand{\B}{\ensuremath{\mathcal{B}}} 
\newcommand{\I}{\ensuremath{\mathcal{I}}}

\newcommand{\vphi}{\ensuremath{\varphi}}    
\newcommand{\vrho}{\ensuremath{\varrho}}
\newcommand{\lalg}[1][g]{\ensuremath{\mathg{#1}}}
\newcommand{\env}[1][g]{\ensuremath{\mathcal{U}(#1)}}
\newcommand{\alg}[1][A]{\ensuremath{\mathcal{#1}}}
\newcommand{\cent}[1][G]{\ensuremath{\mathcal{Z}(#1)}}
\newcommand{\Hom}{\ensuremath{\mathrm{Hom}}}
\newcommand{\End}{\ensuremath{\mathrm{End}}}

\newcommand{\im}{\ensuremath{\mathrm{im}}}

\newcommand{\Lin}{\ensuremath{\mathrm{Lin}}}
\newcommand{\gr}{\ensuremath{\mathrm{gr}\,}}
\newcommand{\Ann}{\ensuremath{\mathrm{Ann}\,}}
\newcommand{\Ind}{\ensuremath{\mathrm{Ind}\,}}
\newcommand{\Res}{\ensuremath{\mathrm{Res}\,}}

\newcommand{\dual}[1][G]{\ensuremath{\widehat{#1}}} 
\newcommand{\Der}{\ensuremath{\mathrm{Der}\,}}

\newcommand{\GL}{\ensuremath{\mathrm{GL}}}

\newcommand{\SO}{\ensuremath{\mathrm{SO}}}

\newcommand{\Diff}[2][]{\ensuremath{\mathcal{D}^{#1}(#2)}}
\newcommand{\Aff}[3][]{\ensuremath{#2^{#1}[#3]}}
\newcommand{\Spec}[1][M]{\ensuremath{\mathcal{S}(#1)}}
\newcommand{\Rat}[3][]{\ensuremath{#2^{#1}(#3)}} 
\newcommand{\Grp}[3][]{\ensuremath{#2^{#1}\lgb #3\rgb }}
\begin{document}
\setcounter{equation}{0}
%
%
\thispagestyle{empty}
%
\date{\today}
\title[Die Frobenius-Zerlegung einer $G$-Variet"at]{Invariante 
Differentialoperatoren und die Frobenius-Zerlegung einer $G$-Variet"at}
%
%
%
\author{Ilka Agricola}
\address{\hspace{-5mm} 
{\normalfont\ttfamily agricola@mathematik.hu-berlin.de}\newline
Institut f\"ur Reine Mathematik \newline
Humboldt-Universit\"at zu Berlin\newline
Sitz: Ziegelstr. 13 A\newline
D-10099 Berlin\\
Deutschland}
\thanks{ Diese Arbeit wurde gef"ordert vom SFB 288 "`Differentialgeometrie 
und Quantenphysik"' der DFG}
\begin{abstract}
Sei $G$ eine zusammenh"angende reduktive komplexe algebraische Gruppe, die auf 
einer glatten affinen komplexen Variet"at $M$ wirke, und bezeichne
$\Diff[G]{M}$ die $G$-invarianten algebraischen Differentialoperatoren auf $M$.
Zerlegt man den Koordinatenring $\Aff{\C}{M}$ in $G$-isotypische Komponenten,
so zeigen wir, da"s die hierbei auftretenden Vielfachheitenr"aume
irreduzible $\Diff[G]{M}$-Moduln sind, zentralen Charakter haben und durch
diesen eindeutig bestimmt sind. Anschlie"send beschreiben wir die analoge
Zerlegung f"ur reelle Formen und zeigen anhand einiger singul"arer Beispiele,
da"s f"ur nicht glatte Variet"aten "ahnliche Ergebnisse nicht zu erwarten
sind.
\end{abstract}
\maketitle
\pagestyle{headings}
%
%
%
\section{Die algebraischen Differentialoperatoren}\label{section-alg-diff-ops}
\noindent
Betrachte eine komplexe reduktive zusammenh"angende algebraische Gruppe $G$,
die auf einer ebenfalls komplexen irreduziblen affinen glatten Variet"at 
$M$ regul"ar wirke. Mittels Translationen operiert $G$ dann auch auf dem 
Koordinatenring $\Aff{\C}{M}$ 
 \bdm
 \vrho(g)f(m)\ = \ f(g^{-1}m) \text{ f"ur }f\in\Aff{\C}{M},\ g\in G\, ,
 \edm
und es ist wohlbekannt, da"s diese $G$-Wirkung
lokal-endlich ist, d.h.\ da"s jeder endlich-dimen\-sion\-ale Unterraum 
von  $\Aff{\C}{M}$ in einem endlich-dimensionalen $G$-invarianten 
Teilraum von $\Aff{\C}{M}$ enthalten ist. Mit diesem Beispiel vor Augen
vereinbaren wir folgende
\begin{dfn}
Mit einem Vektorraum $L$ meinen wir immer einen $\C$-Vektorraum abz"ahlbarer
Dimension. Ein solcher wird \emph{halbeinfacher $G$-Modul} genannt, wenn er
unter $G$ lokal-endlich und vollst"andig reduzibel ist.
\end{dfn}\noindent
Mit $L$ tr"agt auch $\End(L)$ eine Wirkung von $G$; das folgende Beispiel
zeigt jedoch, da"s die Eigenschaft, halbeinfacher $G$-Modul zu sein, sich 
i.a.\ nicht von $L$ auf $\End(L)$ "ubertr"agt. Wir beschr"anken uns
deswegen auf solche Unteralgebren $\alg$ von $\End(L)$, die unter der 
induzierten $G$-Wirkung halbeinfach im eben genannten Sinne sind, und 
werden in K"urze sehen, da"s im Fall $L=\Aff{\C}{M}$ die
algebraischen Differentialoperatoren $\Diff{M}$ diese Eigenschaft haben.
\begin{exa}
Wir w"ahlen $L=\Aff{\C}{x}$ den Polynomring in einer Variablen mit der 
multiplikativen Wirkung von $\C^*$ im Argument. Auf $L$ wird durch
 \bdm
 1 \lmapsto 1,\quad x\lmapsto x^2,\quad x^2\lmapsto x^3+x^4,\ldots,\quad
 x^n\lmapsto x^{n+1}+\ldots + x^{2n}
 \edm
eine lineare Abbildung $T$ definiert, und $x^n$ transformiert sich unter
$T_{\lambda}:= \vrho(\lambda)T\vrho(\lambda)^{-1}$ wie
 \bdm
 T_{\lambda} x^n\ =\ 
 \lambda x^{n+1}+\lambda^2 x^{n+2}+\ldots +\lambda^n x^{2n},\quad n\neq 0\,.
 \edm
Sei $\B$ die von den $T_{\lambda}$, $\lambda\in\C^*$, erzeugte Unteralgebra 
von $\End(L)$. Nach Konstruktion tr"agt sie die von $ \Aff{\C}{x}$
induzierte $\C^*$-Wirkung, f"ur die 
$\vrho(\lambda)T_{\mu}\vrho(\lambda)^{-1}=T_{\lambda\mu}$ gilt.
Leicht sieht man ein, da"s die $T_{\lambda}$ linear unabh"angig sind;
insbesondere ist $\B$ unter der Wirkung von $\C^*$ nicht lokal-endlich. 
\end{exa}\noindent
Zu den "`algebraischen"' Differentialoperatoren einer affinen
Variet"at $M$ gibt es mehrere verschiedene Zug"ange: "uber den Grad als
Differentialoperator, "uber die algebraischen Differentialoperatoren 
$\Diff{V}$ des umgebenden Vektorraums $V$, welche man einfach als
die Differentialoperatoren mit polynomialen Koeffizienten definiert,
"uber die Derivationen von $M$ oder aber "uber die Symbolabbildung.
Diese Definitionen gelten meist -- entsprechend angepa"st -- auch f"ur
Variet"aten mit Singularit"aten, doch liefern sie nur im glatten Fall alle
dieselbe Algebra. Einige dieser verschiedenen Bilder wollen wir kurz 
skizzieren; 
f"ur Beweise sei auf die ausf"uhrliche Beschreibung in dem Lehrbuch
\cite{McConnell&R} verwiesen.

Sei zun"achst $M$ eine nicht notwendigerweise glatte, aber
irreduzible reelle oder komplexe affine Variet"at, und $\alg$ entweder
ihr Koordinatenring $\Aff{\K}{M}$ oder dessen Quotientenk"orper
$\Rat{\K}{M}$. Wir definieren die algebraischen
Differentialoperatoren von $\alg$ als diejenige  filtrierte Algebra,
die als die Vereinigung "uber $k$ aller "`Endomorphismen vom Grad $\leq k$"'
entsteht, also
 \bea[*]
 \Diff{\alg}_0 & := &  \{ T\in \Hom(\alg) \ |\ \qqs f\in \alg:
 \ T f-f T = 0\, \} \, ,\\
 \Diff{\alg}_k & := & \{ T\in \Hom(\alg) \ |\ \qqs f\in \alg:
 \ T f-f  T \in\, \Diff{M}_{k-1} \}\, ,
 \eea[*]
und
 \begin{equation}\label{diff-def}
 \Diff{\alg}\ :=\ \bigcup_{k\geq 0}\Diff{\alg}_k\, .
 \end{equation}
Insbesondere sind die Endomorphismen vom Grad Null genau die
Elemente von $\alg$ und diejenigen vom Grad Eins deren
Derivationen  (vgl.\ \cite[Lemma 15.5.3]{McConnell&R})
 \bdm
 \Diff{\alg}_0\ =\ \alg,\quad \Der \alg\  =\ 
 \{D\in\Diff{\alg}_1\ |\ D(1)=0 \}\, .
 \edm
Dabei verstehen wir ganz allgemein unter den Derivationen $\Der \alg$
diejenigen $\K$-linearen Abbildungen von $\alg$ nach $\alg$, die die 
Leibniz-Regel erf"ullen. 
F"ur uns von Interesse ist die Klasse von
Operatoren, die f"ur $\alg =\Aff{\K}{M}$ ensteht. Wir werden diese
kurzerhand die  \emph{algebraischen Differentialoperatoren von $M$}
nennen und vereinbaren die Bezeichnungen
 \bdm
 \Diff{M}\ :=\ \Diff{\Aff{\K}{M}},\quad \Der M\ :=\ \Der \Aff{\K}{M}\, .
 \edm
Der Fall $\alg=\Rat{\K}{M}$ wurde hier nur deswegen 
mitdefiniert, weil wir ihn zur Formulierung des nun folgenden Kriteriums 
brauchen, mit dem man entscheiden kann, wann ein vorgegebener 
Differentialoperator  algebraisch ist.
\begin{lem}[Erste Charakterisierung von $\Diff{M}$ 
{\cite[Thms.\ 15.1.24, 15.5.5]{McConnell&R}}]\label{diff-ops-krit}
Sei $M$ eine affine irreduzible Variet"at. Dann gilt:
\begin{enumerate}
\item $\Der M = \{D \in \Der \Rat{\K}{M} \ |\ D(\Aff{\K}{M})\subset
 \Aff{\K}{M}\}$;
\item $\Diff{M} = \{D\in\Diff{\Rat{\K}{M}}\ |\ D(\Aff{\K}{M})\subset
 \Aff{\K}{M}\}$.
\end{enumerate}
\qed
\end{lem}
%
\noindent
Sei andererseits $\Diff{M}^*$ die von $\Aff{\K}{M}$ und allen ihren
Derivationen $\Der M$ erzeugte Algebra, filtriert nach der Anzahl Derivationen,
deren Produkt man nimmt. Sie ist immer eine Unteralgebra von
$\Diff{M}$, und ihre nat"urliche Graduierung stimmt mit der von 
$\Diff{M}$ induzierten "uberein.
Weiterhin induziert die faserweise Skalarmultiplikation eine Graduierung
auf dem Koordinatenring $\Aff{\K}{T^* M}=\bigop_n\Aff{\K}{T^* M}^n$  des
Kotangentialb"undels von $M$.
Wie dies aus der Theorie der Differentialoperatoren "uber Mannigfaltigkeiten
zu erwarten ist, stellt das Hauptsymbol eine Bijektion zwischen
den genannten graduierten Ringen her, welche wir nun 
definieren wollen. Sei hierzu $(m,\xi)\in T^* M$, und $g\in\Aff{\K}{M}$   
eine Funktion mit $dg(m)=\xi$. Dann ist das
Hauptsymbol des algebraischen Differentialoperators $D$ vom Grad $k$
gegeben durch
 \bdm
 \sigma_k(D)(m,\xi)\ =\ D\left((g-g(m))^k/k!\right)(m)\, .
 \edm
%
\begin{lem}[Zweite Charakterisierung von $\Diff{M}$]
\label{diffop-isom-kotangbundel}
Sei $M$ eine glatte irreduzible affine Variet"at. Dann ist die
Symbolabbildung ein Algebrenisomorphismus
 \bdm
 \gr \Diff{M}^*\ \cong\ \Aff{\K}{T^* M}
 \edm
und es gilt $\Diff{M}=\Diff{M}^*$.
\end{lem}
\begin{proof}
Dies folgt letztlich aus der Exaktheit der Sequenz (vgl.\ 
\cite[Prop.\ 15.4.9]{McConnell&R})
 \bdm
 0\ \lra\ \Diff{M}_{k-1}\lra\ \Diff{M}_{k}\stackrel{\sigma_k}{\lra}\
 \Aff{\K}{T^* M}^k\lra\ 0\, ,
 \edm
deren Beweis auf der Existenz einer geeigneten Globalisierung beruht.
Mit einem analogen Argument beweist man die zweite Behauptung
\cite[Cor.15.5.6.]{McConnell&R}.
\end{proof}
\begin{NB}
Falls die Variet"at $M$ singul"are Punkte hat, so stimmen die
Algebren $\Diff{M}$ und $\Diff{M}^*$ im allgemeinen nicht "uberein, und 
es wird au"serordentlich schwierig, "uber $\Diff{M}$ irgendwelche 
Aussagen zu machen. Beispiele hierzu werden wir ausf"uhrlich bei der
Diskussion der Frobenius-Zerlegung f"ur singul"are Variet"aten behandeln. 
\end{NB}
Den ersten Schritt f"ur das Studium der $G$-Wirkung auf den
Differentialoperatoren von $M$ bildet die nun folgende Proposition, die
unabh"angig von der Wahl des Grundk"orpers gilt.
\begin{prop}[Einfachheit von $\Diff{M}^*$ {\cite[Thm.\ 15.3.8]{McConnell&R}}]
\label{diff-simple-alg}
Folgende Bedingungen sind \linebreak "aquivalent:
\begin{enumerate}
\item Die Variet"at $M$ ist glatt;
\item Die Algebra $\Diff{M}^*$ ist einfach;
\item \Aff{\K}{M} ist ein irreduzibler $\Diff{M}^*$-Modul.
\end{enumerate}\qed
\end{prop}
\noindent
Ist nun $M$ im Fall $\K=\C$ eine $G$-Variet"at, so wird die
Wirkung von $G$ auf $\Diff{M}$ definiert via 
 \begin{equation}\label{dfn-group-action-on-Diff}
 g\cdot T \ =\ \vrho(g) T \vrho(g)^{-1}\ 
 \text{ f"ur }g\in G,\ T\in \Diff{M},
 \end{equation}
wobei mit $\vrho$ die $G$-Translation auf $\Aff{\C}{M}$ bezeichnet. 
Bez"uglich dieser Wirkung ist die Symbolabbildung $G$-"aquivariant;
somit ist die Wirkung von $G$ auf $\Diff{M}$ wieder lokal $G$-endlich.
Die $G$-invarianten Differentialoperatoren werden mit $\Diff[G]{M}$
bezeichnet und  enthalten die $G$-invarianten Funktionen $\Aff{\C}{M}^G$.
Weil sowohl die Filtrierung nach dem Grad von $\Diff{M}$ als auch die
Graduierung von $\Aff{\C}{T^* M}$ $G$-stabil sind und $\Aff{\C}{T^*M}^G$
unter $G$ lokal-endlich ist, 
gilt das Analogon von Lemma~\ref{diffop-isom-kotangbundel}:
\begin{lem}[Charakterisierung von ${\Diff[G]{M}}$]
\label{invdiffop-isom-invkotangbundel}
Sei $M$ eine komplexe glatte irreduzible affine Variet"at, auf der die
komplexe reduktive zusammenh"angende algebraischen Gruppe $G$ wirke. Dann ist 
 \bdm
 \gr \Diff[G]{M}\ \cong\  \Aff{\C}{T^* M}^G\, .
 \edm
\qed
\end{lem}
%
\begin{exa}\label{exa-env-alg}
Ist $M=G$ mit der linksregul"aren Wirkung von $G$ auf sich selbst und
identifiziert man wie "ublich $\lalg[g]$ mit den linksinvarianten
Vektorfeldern, so ist $T^* G\cong G\x \lalg^*$, also 
$\Aff{\C}{T^* G}\cong\Aff{\C}{G}\ox\Aff{\C}{\lalg^*}$.
Da es auf $G$ keine nichttrivialen $G$-invarianten Funktionen gibt, 
andererseits alle Elemente aus der symmetrischen Algebra $\Aff{\C}{\lalg^*}$
von $\lalg$ unter $G$ invariant sind, ist 
$\gr \Diff[G]{G}\cong\Aff{\C}{\lalg^*}$;
bekannterweise ist dies als graduierte Algebra isomorph zur universellen 
Einh"ullenden $\env[\lalg]$, die uns als Algebra der
$G$-invarianten Differentialoperatoren auf $G$ wohlvertraut ist.

Andererseits weicht die Struktur von $\Diff{M}$ (oder $\Diff[G]{M}$) im
Detail doch sehr von der
einer universellen Einh"ullenden ab; denn wenn $\Diff{M}$
(bzw.\ $\Diff[G]{M}$) als filtrierte Algebra zur einh"ullenden Algebra einer 
Lie-Agebra $\lalg[k]$ isomorph ist,
dann mu"s $\Aff{\C}{T^*M}$ (resp.\ $\Aff{\C}{T^*M}^G$) deren
symmetrische Algebra $\Aff{\C}{\lalg[k]^*}$ sein, also auf jeden Fall ein
Polynomring. Dies ist im allgemeinen sicher nicht der Fall, wenn $M$ kein
Vektorraum oder $T^*M$ kein triviales B"undel ist. 
\end{exa}\noindent
Beim Studium der Unteralgebren der Endomorphismen
halbeinfacher $G$-Moduln ben"otigt man immer
wieder folgende, Dixmier zugeschriebene Variante des Schurschen Lemmas f"ur
Vektorr"aume abz"ahlbarer Dimension. Vorher allerdings noch eine
Vereinbarung:
\begin{dfn}
Es sei fortan und bis zum Ende von Abschnitt \ref{section-central-char}
der Grundk"orper gleich $\C$, auch wenn dies nicht immer erw"ahnt wird.
\end{dfn}
\begin{lem}[Lemma von Dixmier]\label{Dixmier-lemma}
Sei $V$ ein $\C$-Vektorraum abz"ahlbarer Dimension.
\begin{enumerate}
\item 
Zu jedem $T\in\End (V)$ existiert ein $q\in\C$ so, da"s $T-qI$ nicht
invertierbar ist;
\item Ist $\alg$ eine Teilalgebra von $\End(V)$, die auf $V$ irreduzibel
wirkt, so ist $\End_{\alg}(V)\, =\, \C$.
\end{enumerate}
\end{lem}
\begin{proof}
Die erste Behauptung findet sich in \cite[S.11]{Wallach2}. 
F"ur die zweite sei $Z\in \End_{\alg}(V)$ und $q\in\C$
derart, da"s $Z-qI$ nicht invertierbar ist, d.h.\ es ist
entweder $\ker (Z-qI)\neq \{0\}$ oder $\im (Z-qI)\neq V$.
Weil aber $Z$ mit allen Elementen aus $\alg$ vertauscht, m"ussen
$\ker (Z-qI)$ und $\im (Z-qI)$ beide $\alg$-invariant sein; aus der
Irreduzibilit"at von $V$ folgt damit, da"s entweder
$\ker (Z-qI)=V$ oder $\im (Z-qI)=\{0\}$ sein mu"s. In beiden F"allen
folgt $Z=qI$ auf ganz $V$.
\end{proof}
\begin{exa}
Sei $V$ halbeinfacher $G$-Modul.
Man sagt, da"s $V$ \emph{zentralen Charakter} hat, wenn jedes Element $z$
im Zentrum $\cent[\lalg[g]]$ der universellen Einh"ullenden von $G$ skalar 
wirkt. Unter Verwendung des 
Harish-Chandra-Homomorphismus $\cent[\lalg[g]]\cong \Aff{\C}{\lalg[h]^*}^W$
ist z.B.\ f"ur die endlich-dimensionale $G$-Darstellung $V^{\lambda}$
vom h"ochsten Gewicht $\lambda$ der zentrale Charakter gleich
$\lambda+\vrho$, wobei $\vrho$ die H"alfte der positiven Wurzeln von $G$
sei.
\end{exa}
\begin{dfn}
Aus dem Lemma von Dixmier folgt nun, da"s
f"ur einen irreduziblen $\alg$-Modul $U$, $\alg$ wieder eine Unteralgebra von
$\End(U)$, das Zentrum von $\alg$ auf $U$ skalar wirkt, genauer,
es existiert ein Homomorphismus $\chi: \cent[\alg]\ra\C$ mit 

 \bdm
 Z\cdot u\ =\ \chi(Z) \cdot u \quad \qqs u\in U\, ,
 \edm
den wir den \emph{zentralen Charakter} von $\alg$ auf $U$ nennen.
%
\end{dfn}
Mit Hilfe dieses Lemmas k"onnen wir das Zusammenspiel der
algebraischen Differentialoperatoren mit den Endomorphismen 
von $\Aff{\C}{M}$, als deren Teilalgebra sie definiert waren,
noch pr"azisieren. Zusammen mit Proposition \ref{diff-simple-alg}
folgt n"amlich:
\begin{prop}[Schursches Lemma]\label{irred-schur}
Ist $M$ eine glatte irreduzible affine Variet"at, so gilt
 \bdm
 \End_{\Diff{M}}(\Aff{\C}{M})\ =\ \C\, .
 \edm
\qed
\end{prop}
\noindent
Intuitiv bedeutet dies, da"s $\Diff{M}$ von $\End (\Aff{\C}{M})$ nicht
allzu sehr abweichen kann. Auch diese Aussage l"a"st sich mit einem
allgemeinem Satz genauer fassen, den wir deswegen zuerst in der f"ur uns 
relevanten Fassung zitieren wollen.
\begin{thm}[Dichtheitssatz von Jacobson {\cite[Thm, Ch.12, 2]{Pierce}}]
\label{thm-jacobson-abstr}
Sei $L$ ein $\C$-Vektorraum und
$\alg$ eine Teilalgebra von $\End(L)$, die auf $L$ irreduzibel
wirkt. Sind $u_1,\ldots,u_n$ linear unabh"angige Elemente  und
$v_1,\ldots,v_n$ beliebige Elemente von $L$, so existiert ein 
$T\in\alg$ derart, da"s $T u_i=v_i$ f"ur alle $i$ gilt. 
Das bedeutet:
ist $X$ ein endlich-dimensionaler Unterraum von $V$, so gilt
 \bdm
 \alg\restr_X\ =\ \Hom(X,L)\, .
 \edm
Insbesondere gilt also mit Proposition~\ref{irred-schur}, da"s 
 \bdm
 \Diff{M}\restr_{X}\ =\ \Hom(X, \Aff{\C}{M})
 \edm
f"ur jede irreduzible glatte affine Variet"at $M$ ist.
\qed
\end{thm}
%
\begin{prop}\label{thm-jacobson}
Sei $G$ eine reduktive algebraische Gruppe, $L$ ein halbeinfacher $G$-Modul,
und $\alg$ eine Unteralgebra von $\End(L)$, die auf $L$ irreduzibel wirkt
und bez"uglich der induzierten $G$-Wirkung auf $\alg$ ebenfalls halbeinfach 
ist.
Sei weiterhin $X$ ein endlich-dimensionaler $G$-invarianter Unterraum von $L$.
Dann gilt f"ur die Algebra $\alg^G$ der $G$-Invarianten von $\alg$:
 \bdm 
 \alg^G\restr_{X}\ =\ \Hom_G (X,L);
 \edm
insbesondere ist also
 \bdm
 \Diff[G]{M}\restr_{X}\ =\ \Hom_G (X, \Aff{\C}{M})
 \edm
f"ur eine irreduzible glatte Variet"at $M$.
%
\end{prop}
\begin{proof}
Sei $T\in \Hom_G(X, L)$. Nach Satz~\ref{thm-jacobson-abstr} existiert 
zun"achst ein $D\in\alg$ mit $D\restr_{X}=T$. Nach Voraussetzung 
liegt $D$ in
einem endlich-dimensionalen $G$-invarianten Unterraum $E$ von $\alg$.
Wegen der Reduktivit"at von $G$ zerf"allt $E$  in $G$-invariante
Teilr"aume $E=E^G \op F$; weil $T$ zudem $G$-invariant war, ist 
$T\in E^G\restr_{X}$, es mu"s also ein $\tilde{D}\in E^G\subset \Diff[G]{M}$
mit $\tilde{D}\restr_{X}=T$ geben. 

F"ur den Fall $L=\Aff{\C}{M}$ und
$\alg=\Diff{M}$ verbleibt deswegen blo"s, die Voraussetzungen des Satzes
zu "uberpr"ufen. Es sind $\Aff{\C}{M}$ und $\Diff{M}$, mittels Lemma
\ref{invdiffop-isom-invkotangbundel} mit $\Aff{\C}{T^* M}$ identifiziert,
lokal-endlich bzgl.\ der Wirkung von $G$.  Proposition~\ref{diff-simple-alg} 
stellt sicher, da"s  $\Diff{M}$ auf $\Aff{\C}{M}$ irreduzibel wirkt.
\end{proof}
%
Abschlie"send stellen wir hier noch einige weitere Eigenschaften der
algebraischen Differentialoperatoren zusammen, die Friedrich Knop
im Laufe der letzten Jahre in einer Reihe von Arbeiten entwickelte. 
Es bezeichne hierbei $\alg'$ immer den Kommutanten einer 
Unteralgebra $\alg$ von $\Diff{M}$. Mit $\env[\lalg[g]]$ ist hier
stillschweigend dessen homomorphes Bild in $\Diff{M}$ gemeint.
\begin{thm}[Satz vom Bikommutanten {\cite[Thm.9.5.]{Knop94}}]
\label{alg-bikommutant}
Die zusammenh"angende reduktive komplexe algebraische Gruppe $G$ wirke
auf der  glatten affinen Variet"at $M$.  Dann gilt:	
\begin{enumerate} 
\item  $\Diff[G]{M}\ =\ \Diff[G]{M}''$, woraus f"ur die Zentren sofort
die Beziehung
 \bdm
 \cent[{\Diff[G]{M}'}]\ =\ \Diff[G]{M}'\cap\Diff[G]{M}\ =\ 
\cent[{\Diff[G]{M}}]\ =: \cent[M]
 \edm
folgt;
\item Es existiert ein zentraler Differentialoperator $0\neq T\in \cent[M]$
derart, da"s die Inklusionen
 \bdm
 T\cdot \Diff[G]{M}'\ \subset\ \env[\lalg[g]] \cent[M]\ \subset\ \Diff[G]{M}'
 \edm
gelten;
\item Die Multiplikation ist ein Isomorphismus
\bdm
\Diff[G]{M}'\ox_{\cent[M]}\Diff[G]{M}\stackrel{\cong}{\lra}\cent[M]',\quad
a\ox b \lmapsto ab=ba;
\edm
\item Die Algebren $\Diff{M}$, $\Diff[G]{M}$ und $\Diff[G]{M}'$ sind freie
$\cent[M]$-Moduln (im ersten Fall sowohl als Links- als auch als Rechtsmodul).
\end{enumerate}\qed
\end{thm}
%
\begin{dfn}
Fortan bezeichne $\cent[M]$ immer das Zentrum von $\Diff[G]{M}$. 
F"ur $M=G$ ist dies genau das Zentrum der universellen Einh"ullenden
von $\lalg[g]$, weswegen wir hierf"ur inkonsequenterweise $\cent[\lalg[g]]$
statt $\cent[G]$ schreiben werden.
\end{dfn}
\begin{NB}
Dieser Satz wird beim Studium der zentralen Charaktere von $\Diff[G]{M}$ 
von Bedeutung sein. Vorerst sei kurz angemerkt, wie die zweite Aussage
zu verstehen ist: sie besagt, da"s  $\Diff[G]{M}'$ und $\env[\lalg[g]]$
nur um Elemente in $\cent[M]$ voneinander abweichen. Insbesondere liefert
sie eine Kontrolle von $\Diff[G]{M}'$, wenn man die Wirkung von $\cent[M]$
kennt.
\end{NB}
\section{Die abstrakte Frobenius-Zerlegung}\noindent
Es bezeichne
$B=HN$ eine Boreluntergruppe ($H$ eine Cartanuntergruppe,  $N$ ein 
unipotentes Radikal von $B$)
und $\dual =\{(V^{\lambda},\, \pi_{\lambda})\}$ die Menge
aller "Aquivalenzklassen endlich-dimensionaler
Darstellungen von $G$, die wir stillschweigend den dominanten 
Gewichten $\lambda\in P_+(G)$ gleichsetzen werden. Ist $\mu\in P(G)$ ein
beliebiges Gewicht von $G$, dann bezeichen wir den zugeh"origen Charakter 
von $H$ mit $h^{\mu}$ und erweitern ihn zu einem Charakter von $B$
mittels $(hn)^{\mu}:=h^{\mu}$ ($h\in H$, $n\in N$).
Ist nun $L$ ein halbeinfacher $G$-Modul, so werden die 
$N$-invarianten Vektoren vom Gewicht $\mu$ definiert als
 \bdm
 L^N(\mu)\ =\ 
 \{f\in L\ |\ \ \vrho(b) f\, =\, b^{\mu} f \quad \qqs\, b\in B\}.
 \edm
Wenn $\mu=\lambda$ dominant ist, sind diese isomorph (als Vektorraum) zu den
$G$-invarianten Morphismen
von $V^{\lambda}$, dem irreduziblen Modul zum h"ochsten Gewicht $\lambda$,
nach $L$
 \bdm
 L^N(\lambda)\ \cong\ \Hom_G (V^{\lambda}, L),
 \edm
wie man durch Auswertung auf einem Vektor vom h"ochsten Gewicht $\lambda$
leicht einsieht. Wir definieren eine Abbildung
 \bdm
 \vphi_{\lambda}:\ \Hom_G (V^{\lambda}, L)\ox V^{\lambda}
 \lmapsto L,\quad \vphi_{\lambda}(T\ox v):= T(v),
 \edm
die die Eigenschaft hat, $G$-"aquivariant zu sein, wenn man
 \bdm
 g\cdot (T\ox v)\ :=\ T\ox\pi_{\lambda}(g) v
 \edm
setzt. F"ur eine Teilalgebra $\alg$ von $\End(L)$ wird die Wirkung von 
$D\in\alg$ auf
$T\in\Hom_G(V^{\lambda}, L)$ mittels Komposition definiert,
$D(T)=D\circ T$, und diese Definition impliziert, da"s 
$\Hom_G( V^{\lambda}, L)$ auf jeden Fall $\alg^G$-invariant
ist; denn erf"ullt $T$ die Bedingung 
$T(\pi_{\lambda}(g) v)=\pi_{\lambda}(g)T(v) $, so folgt unmittelbar
f"ur ein $D\in\alg^G$
 \bdm
 (D\circ L)(\pi_{\lambda}(g) v)\ =\ D(L(\pi_{\lambda}(g) v))\ =\
 \pi_{\lambda}(g)D(L(v))\ =\ \pi_{\lambda}(g)(D\circ L)(v) \, .
 \edm
Die Menge
 \bdm
 \Spec[L]\ :=\ \{\lambda\in P_+(G)\ |\ L^N(\lambda)\neq 0\},
 \edm
wird das \emph{Spektrum} von $G$ auf $L$ genannt. Ist nun 
$L=\Aff{\C}{M}$ der Koordinatenring einer irreduziblen affinen $G$-Variet"at,
so schreibt man der K"urze halber $\Spec$ statt
$\Spec[\Aff{\C}{M}]$ und kann hier"uber noch folgende Aussagen machen:
weil unter der punktweisen Multiplikation von Funktionen
 \bdm
 \Aff{\C}{M}^N(\lambda)\cdot\Aff{\C}{M}^N(\mu)\, \subset\,
 \Aff{\C}{M}^N(\lambda+\mu)
 \edm
gilt, ist $\Spec$ eine additive Halbgruppe, und mittels des Satzes von 
Rosenlicht zeigt man, da"s sie  als solche immer
endlich erzeugt ist (vgl.\ z.B.\ \cite[Cor.6]{Panyushev90} sowie
die Anmerkungen in Definition~\ref{rank-compl}).
Weil $\Aff{\C}{M}^N(0)=\Aff{\C}{M}^G$ gilt und jeder der R"aume
$\Aff{\C}{M}^N(\lambda)$ "uber $\Aff{\C}{M}^G$ endlich erzeugt ist, sind 
erstere genau dann endlich-dimensional, wenn $\Aff{\C}{M}^G=\C$ ist,
es auf $M$ neben den konstanten also keine weiteren $G$-invarianten 
Funktionen gibt (vgl.\  \cite{Kraft}).
%
\begin{NB}
Der folgende Satz ist eine Verallgemeinerung eines Resultats von N.~Wallach
\cite[Prop.1.5]{Wallach93a},
welches den Fall eines Vektorraums behandelt. Obgleich wir uns von dem 
dortigen Beweis sowie seiner Darstellung in \cite[Thm.4.5.12]{Goodman&W}
haben inspirieren lassen, liegt der wesentliche Unterschied hier
in der Tatsache, da"s die Graduierung der Polynome und Differentialoperatoren
mit polynomialen Koeffizienten nicht mehr gebraucht wird. 
Eine Ank"undigung dieses Ergebnisses ist zudem implizit und ohne Beweis
in \cite[S.219]{Bien93} enthalten. 
\end{NB}

\begin{thm}[Abstrakte Frobenius-Zerlegung]\label{isotyp-zerl-abstr}
Sei $G$ eine zusammenh"angende komplexe reduktive algebraische 
Gruppe, $L$ ein halbeinfacher $G$-Modul,
und $\alg$ eine Unteralgebra von $\End(L)$, die auf $L$ irreduzibel wirkt
und bez"uglich der induzierten $G$-Wirkung ebenfalls halbeinfach ist.
Dann zerf"allt $L$ als $G$-Modul in
 \bdm
  L\ \cong\ 
 \bigop_{\lambda\in P_+(G)} \Hom_G(V^{\lambda}, L)\ox V^{\lambda},
 \edm
und es gilt:
\begin{enumerate}
\item $\Hom_G(V^{\lambda}, L)\cong L^N(\lambda)$ ist ein irreduzibler
$\alg^G$-Modul;
\item F"ur $\lambda\neq\mu$ sind die $\Hom_G(V^{\lambda}, L)$
als $\alg^G$-Moduln paarweise nicht isomorph;
\item Die R"aume $\Hom_G(V^{\lambda}, L)\cong L^N(\lambda)$ haben,
aufgefa"st als $\alg^G$-Moduln, zentralen Charakter. 
\end{enumerate}
\end{thm}
%
\begin{kor}[Frobenius-Zerlegung von $\Aff{\C}{M}$]\label{isotyp-zerl-aff}
Es sei $G$ eine zusammenh"angende komplexe algebraische reduktive Gruppe, die
auf der glatten affinen irreduziblen Variet"at $M$ wirkt. Dann zerlegt sich
der Koordinatenring von $M$ als $G$-Modul in
 \bdm
 \Aff{\C}{M}\ \cong\ 
 \bigop_{\lambda\in P_+(G)} \Hom_G(V^{\lambda}, \Aff{\C}{M})\ox V^{\lambda},
 \edm
und es gilt:
\begin{enumerate}
\item $\Hom_G(V^{\lambda}, \Aff{\C}{M})\cong\Aff{\C}{M}^N(\lambda)$ ist 
ein irreduzibler $\Diff[G]{M}$-Modul;
\item F"ur $\lambda\neq\mu$ sind die $\Hom_G(V^{\lambda}, \Aff{\C}{M})$
als $\Diff[G]{M}$-Moduln paarweise nicht isomorph;
\item Die R"aume $\Hom_G(V^{\lambda}, \Aff{\C}{M})\cong\Aff{\C}{M}^N(\lambda)$
haben, aufgefa"st als $\Diff[G]{M}$-Moduln, zentralen Charakter. 
\end{enumerate}
\end{kor}
\begin{proof}
Die Zerlegung an sich ist ein Standard-Ergebnis und kann z.B.\ in
\cite[Thm.12.1.3]{Goodman&W} nachgelesen werden. 
F"ur die eigentlichen Behauptungen verfahren  wir  wie folgt:
\begin{enumerate}
\item 
Sei $X^{\lambda}$ ein zu $V^{\lambda}$ isomorpher Unterraum von 
$L$ sowie $f\neq 0$ ein nicht-triviales Element aus $X^{\lambda}$.
Wir studieren den $\alg^G$-invarianter Teilraum 
$\alg^G\cdot f=: U_f$.
Dann ist nat"urlich $U_f\cap X^{\lambda}\neq \{0\}$, denn alle Vielfachen
von $f$ liegen sowohl in der linken als auch in der rechten Seite.
Dies ist bereits der ganze Durchschnitt: ist $u\in U_f\cap X^{\lambda}$
beliebig und $D$ ein Operator aus $\alg^G$ derart, da"s
$Df=u$ ist, so folgt aus der $G$-Irreduzibilit"at von $X^{\lambda}$, da"s
$X^{\lambda}$  sowohl gleich $\Lin(G\cdot u)$ als auch gleich $\Lin(G\cdot f)$
ist. Deswegen existieren f"ur jedes $x\in X^{\lambda}$ Gruppenelemente 
$g_i\in G$ derart, da"s $x=\sum\alpha_i\vrho(g_i)f$ ist. Also ist
 \bdm
 Dx\, =\, \sum \alpha_i D\vrho(g_i)f \, =\, \sum \alpha_i \vrho(g_i)Df \, =
 \sum \alpha_i \vrho(g_i)u,
 \edm
es ist somit $D$ ein $G$-invarianter Homomorphismus von $X^{\lambda}$
nach $X^{\lambda}$, nach dem Schurschen Lemma demnach ein Vielfaches
der Identit"at. Damit ist $U_f\cap X^{\lambda}=\C\cdot f$ wie behauptet.

Es ist $U_f$ als $\alg^G$-Modul irreduzibel. In der Tat, sei $p\neq 0$ ein 
Element eines nicht-trivialen $\alg^G$-invarianten Unterraums $E$ von
$U_f$. Aus der lokalen $G$-Endlichkeit von $L$ folgt, da"s $p$ in einem
endlich-dimensionalen $G$-invarianten Unterraum, der  $X$ hei"sen m"oge, 
liegen mu"s. Als $G$-Darstellung mu"s $X$ isomorph zu einem Vielfachen von 
$V^{\lambda}$ sein, d.h.\ es existiert ein $G$-Homomorphismus 
$T:X\lra X^{\lambda}$ mit $Tp\neq 0$. Zudem wird nach dem Dichtheitssatz
von Jacobson (Proposition~\ref{thm-jacobson}) $T$ durch einen Operator
$D\in\alg^G$ dargestellt; weil $p$ aber in einem $\alg^G$-invarianten
Unteraum lag, ist $Dp\in E$, und deswegen $E\cap X^{\lambda}\neq \emptyset$. 
Aus $U_f\cap X^{\lambda}=\C\cdot f$ folgt dann, da"s $f$ ebenfalls in $E$
liegt, also ist $U_f$ wie behauptet irreduzibel.

Nun sei $f_1,\ldots,f_n$ eine Basis von $X^{\lambda}$. Die $U_{f_i}$,
die wir fortan mit $U_i$ abk"urzen wollen,
sind  als $\alg^G$-Moduln paarweise isomorph: Nach dem Satz von 
Burnside gibt es Elemente $u_{ji}\in \Grp{\C}{G}$, die $f_i$ auf $f_j$
abbilden, $\vrho(u_{ji})f_i=f_j$, also ist
 \bdm 
 u_{ji}:\quad U_{i}\lra U_{j}
 \edm
ein $\alg^G$-Modul-Isomorphismus. Sodann zeigen wir, da"s die Summe
$\sum U_i$ direkt ist. Wieder nach dem Satz von Burnside ist jeder
Endomorphismus von $X^{\lambda}$ das Bild unter $\vrho$ eines Elements
aus $\Grp{\C}{G}$, der Gruppenalgebra von $G$, wir finden also
Elemente $u_i\in\Grp{\C}{G}$, die $\vrho(u_i)f_j=\delta_{ij}f_i$ erf"ullen. 
Seien nun
die $m_i$ jeweils in $U_i$ mit $\sum m_i=0$, und weiterhin $D_i\in\alg^G$ 
Operatoren mit $D_if_i=m_i$. Dann folgt insgesamt
 \bdm
 0\, =\, \vrho(u_j)(\sum_{i}m_i)\, =\, \sum_{i}\vrho(u_j)D_if_i\, =\,
 \sum_{i}D_i\vrho(u_j)f_i\, =\, D_j f_j\, =\, m_j \quad \qqs j=1,\ldots,n,
 \edm
und die Summe ist  direkt.

Ist nun $Z^{\lambda}$ ein weiterer zu $V^{\lambda}$ isomorpher Unterraum
von $L$, so liegt auch $Z^{\lambda}$ in $\bigop_{i} U_i$: Weil die $f_i$
eine Basis von $X^{\lambda}$ bilden, existiert ein $G$-Homomorphismus
$T:\ X^{\lambda}\ra Z^{\lambda}$, der die $f_i$ auf eine Basis $Tf_i$
von $Z^{\lambda}$ abbildet, und weil $X^{\lambda}$ endlich-dimensional war,
ist $T$ die Restriktion eines Operators $D$ aus $\alg^G$. Damit ist
$Tf_i\in U_i\cap Z^{\lambda}$ f"ur alle $i$, d.h.\ die Basis
$Tf_1,\ldots Tf_n$ von $Z^{\lambda}$ liegt in $\bigop_{i} U_i$, was behauptet
wurde. Dies bedeutet aber nicht anderes, als da"s  $\bigop_{i} U_i$
die $V^{\lambda}$-isotypische Komponente von $L$, also gleich
$\Hom_G(V^{\lambda}, L)\ox V^{\lambda}$ ist, womit die Irreduzibilit"at
das ersten Faktors unter der $\alg^G$-Wirkung bewiesen ist.
%
\item 
Es bezeichne $L_{\lambda}$ bzw.\ $L_{\mu}$ die vollst"andige isotypische
Komponente von $L$ zum Gewicht $\lambda$ bzw.\ $\mu\in\Spec[L]$,
und $U_{\lambda}$ bzw.\ $U_{\mu}$ einen zu $L^N(\lambda)$ bzw.\ 
$L^N(\mu)$ isomorphen $\alg^G$-Modul darin.
Angenommen, es existiert  ein $\alg^G$-Modul-Isomorphismus
 \bdm
 Q:\ U_{\lambda}\lra U_{\mu}\, .
 \edm
Nach 1.\ k"onnen wir ein $f\in L_{\lambda}$ w"ahlen,
welches $U_{\lambda}$ erzeugt, also
$\alg^G\cdot f= U_{\lambda}$. Wir betrachten den
minimalen $G$-invarianten Teilraum $E$ von $L$, der $f$ und $Qf$
enth"alt. Die Projektion $\pi: E\mapsto L_{\lambda}$ vertauscht mit $G$,
ist also ein Element von $\Hom_G(E, L)$; damit existiert ein
$G$-invarianter Operator $D\in\alg^G$, der $\pi$ auf $E$ darstellt, 
$D\restr_E=\pi$, und es ergibt sich
 \bdm
 \pi Qf\ =\ DQf\ =\ QDf\ =\ Q\pi f\ =\ Qf\, ,
 \edm
es ist also $Qf\in L_{\lambda}$, und damit $V^{\lambda}=V^{\mu}$, was mit der
Gleichheit der Charaktere $\lambda$ und $\mu$ gleichbedeutend ist.
\item 
Die Existenz eines zentralen Charakters f"ur die $\alg^G$-Moduln
$\Hom_G(V^{\lambda},L)$, folgt nun sofort aus 1.\ und dem Lemma von Dixmier
(Lemma~\ref{Dixmier-lemma}), da sie 
-- als $\alg^G$-invariante Unterr"aume von $L$ -- von abz"ahlbarer 
Dimension sind.
\end{enumerate}
\end{proof}
Ein intensiv studierter Spezialfall ist der, wo alle Vielfachheiten gleich
Eins sind ($M$ wird dann \emph{sph"arisch} genannt), und der folgende
Satz von Vinberg und Kimel'fel'd liefert ein elegantes
Kriterien um zu testen, ob dieser Fall vorliegt:
\begin{thm}[Vielfachheitenfreie Wirkungen]
Ist $G$ eine zusammenh"angende komplexe algebraische reduktive Gruppe mit 
Boreluntergruppe $B=HN$ und $M$
eine affine $G$-Variet"at, so sind folgende Bedingungen "aquivalent:
\begin{enumerate}
\item $\Aff{\C}{M}$ ist vielfachheitenfrei, d.h.\
$\dim \Aff{\C}{M}^N(\lambda)=1 \ \text{ f"ur alle } \lambda\in\Spec$;
\item $B$ hat einen offenen Orbit in $M$;
\item Jede rationale $B$-invariante Funktion auf $M$ ist konstant: 
$\Rat{\C}{M}^B=\C$.
\end{enumerate}
\end{thm}
%
\begin{proof}
Die "Aquivalenz von 1.\ und 2.\ ist Satz 2 aus der klassischen Arbeit von
Vinberg-Kimel'fel'd \cite{Vinberg&K78}; man kann die Gleichwertigkeit der 
drei Bedingungen aber auch in \cite[Satz 1, III.3.6, S.199]{Kraft} nachlesen.
\end{proof}
%
\begin{prop}[Glatte vielfachheitenfreie Wirkung und Differentialoperatoren]
\label{multfree-diffops}
Die Wirkung einer zusammenh"angenden komplexen algebraischen reduktiven 
Gruppe $G$ auf einer glatten affinen $G$-Variet"at $M$ ist genau dann
vielfachheitenfrei, wenn die Algebra der $G$-invarianten 
Differentialoperatoren $\Diff[G]{M}$ kommutativ ist.
\end{prop}
Dies folgt sofort aus der Tatsache, 
da"s der Vielfachheitenraum ein irreduzibler $\Diff[G]{M}$-Modul ist: ist
$\Diff[G]{M}$ abelsch, dann kann es nur eindimensionale irreduzible
Darstellungen geben. Sind andererseits alle Vielfachheiten gleichs Eins, dann
mu"s ein beliebiges $D\in\Diff[G]{M}$  mittels Skalarmultiplikation auf den
einzelnen Summanden der Zerlegung wirken, und zwei solche kommutieren
nat"urlich.\qed
\begin{NB}
Dieser Sachverhalt war bereits in vielen Spezialf"allen bekannt, etwa f"ur
Vektorr"aume (vgl.\ \cite{Howe&U91a}, wo $\Diff[G]{M}$ f"ur 
vielfachheitenfreie irreduzible Darstellungen ausgearbeitet ist) oder f"ur 
homogene R"aume. Im singul"aren Fall h"angt die
Antwort davon ab, ob man $\Diff{M}$ oder $\Diff{M}^*$ betrachtet; allgemeine
Ausssagen lassen sich hier nicht beweisen, doch stellen wir in 
Abschnitt \ref{sing-exa} einige Beispiele vor, die beweisen, da"s die
vorangegangene Behauptung f"ur $\Diff{M}^*$ falsch ist. 
\end{NB}
\begin{exa}
Wir betrachten den klassichen Fall $M=G$ mit der linksregul"aren Wirkung, auf 
den sich der aller homogenen 
$G$-R"aume reduzieren l"a"st. Sei $V^{\pi}$ Darstellung von $G$
und $V^{\mu}$ Darstellung einer Untergruppe $H$.
Aus der Frobenius-Reziprozit"at f"ur die induzierte bzw.\ eingeschr"ankte
Darstellung von $\mu$ und $\pi$
 \bdm
 \Hom_G(V^{\pi}, \Ind^G_H (\mu))\ =\ \Hom_H(\Res^G_H(\pi), V^{\mu})
 \edm
folgt f"ur die triviale Darstellung $\mu=\id$ von $H=\{e\}$
 \bdm
 \Hom_G(V^{\pi}, \Aff{\C}{G})\ =\
 \Hom (V^{\pi},\C)\ =\ (V^{\pi})^*,
 \edm
welches nat"urlich f"ur alle $\pi\in P_+(G)$ nicht der Nullraum ist,
es ist also $\Spec[G]=P_+(G)$. Bezeichnet man die positiven Wurzeln von $G$
mit $\Phi^+$, so ist wohlbekannt, da"s die Weylgruppe von $G$ ein
eindeutiges Element $w_0$ mit der Eigenschaft $w_0 \Phi^+=-\Phi^+$ enth"alt,
und da"s das h"ochste Gewicht des zu $V^{\lambda}$ dualen $G$-Moduls
gleich $\lambda^*:=-w_0 \lambda$ ist. Man erh"alt also
 \bdm 
 \Aff{\C}{G}\ \cong\ 
 \bigop_{\lambda\in P_+(G)} V^{\lambda^*}\ox V^{\lambda}\, ,
 \edm
wobei die $V^{\lambda^*}$ irreduzible $\env[\lalg[g]]$-Moduln mit
zentralem Charakter $\lambda^* +\vrho$ sind. 
\end{exa}
\section{Die zentralen Charaktere}\label{section-central-char}\noindent
Nach dem Satz von Racah hat das Zentrum der universellen Einh"ullenden
von $G$ genau $\rk G$ erzeugende Elemente, und jede endlich-dimensionale
$G$-Darstellung ist durch die Eigenwerte dieser Operatoren, der sog.\
Casimir-Operatoren von $G$, eindeutig bestimmt 
(vgl.\ \cite[Thm.7.4.7]{Dixmier2}). Bemerkenswerterweise haben auch
die $\Diff[G]{M}$-Moduln $U_{\lambda}$ aus Korollar \ref{isotyp-zerl-aff} diese
Eigenschaft, und dies unabh"angig davon, ob sie endlich- oder
unendlich-dimensional sind. 
%
\begin{thm}[Eindeutigkeit der zentralen Charaktere]\label{central-char-uniq}
Die Vielfachheitenr"aume $U_{\lambda}=\Aff{\C}{M}^N (\lambda)$ sind,
als $\Diff[G]{M}$-Moduln, eindeutig durch ihren zentralen Charakter
bestimmt.
\end{thm}
\begin{proof}
Der Beweis st"utzt sich auf die isotypische Zerlegung von $\Aff{\C}{M}$
als $\Diff[G]{M}$-Modul (Korollar \ref{isotyp-zerl-aff}) und den
Knopschen Bikommutantensatz (Satz \ref{alg-bikommutant}). 
Sei $T$ ein Element aus dem Kommutanten $\Diff[G]{M}'$ von
$\Diff[G]{M}$. Nach Korollar \ref{isotyp-zerl-aff} mu"s $T$ die
Zerlegung von $\Aff{\C}{M}$ invariant lassen; weil es zudem
mit ganz $\Diff[G]{M}$ vertauscht, wirkt $T$ nur auf die 
$\Diff[G]{M}$-Vielfachheitenr"aume, welche genau die  $G$-Darstellungen
$V^{\lambda}$ sind. Die $V^{\lambda}$ sind als $\Diff[G]{M}'$-Moduln
irreduzibel und paarweise nicht "aquivalent. In der Tat, alle zentralen 
invarianten Differentialoperatoren wirken auf $V^{\lambda}$ skalar. Also 
impliziert Satz \ref{alg-bikommutant} (2),
da"s sich $\Diff[G]{M}'$ und $\env[\lalg[g]]$ bei Wirkung hierauf nicht
unterscheiden:
 \bdm
 \Diff[G]{M}'\restr_{V^{\lambda}}\ =\ \env[\lalg[g]]\restr_{V^{\lambda}}\, .
 \edm
Weil die $V^{\lambda}$ irreduzible und paarweise in"aquivalente
$\env[\lalg[g]]$-Moduln sind, sind sie dies auch als $\Diff[G]{M}'$-Moduln.
Insgesamt ist nun die Frobenius-Zerlegung aus Korollar \ref{isotyp-zerl-aff}
vielfachheitenfrei bzgl.\ der $\Diff[G]{M}\ox\Diff[G]{M}'$-Wirkung.
Nach Satz \ref{alg-bikommutant} (3) ist die Multiplikation aber ein
Isomorphismus zwischen $\Diff[G]{M}\ox\Diff[G]{M}'$ und  dem Kommutanten des
Zentrums $\cent[M]$ der invarianten Operatoren; die Zerlegung ist also
vielfachheitenfrei unter der Wirkung der (sehr gro"sen) Algebra $\cent[M]'$,
was zur Folge hat, da"s sich die zentralen Charaktere paarweise unterscheiden
m"ussen.
\end{proof}
Wir erg"anzen nun die Ergebnisse von Friedrich Knop zur Struktur der
$G$-invarianten Differentialoperatoren. Um diese
formulieren zu k"onnen, erinnern wir an dieser Stelle an die Definition 
des Rangs und der Komplexit"at einer affinen $G$-Variet"at.
\begin{dfn}\label{rank-compl}
Den Arbeiten von E.\ Vinberg und D.\ Panyushev folgend, 
definieren wir (f"ur einen Punkt $x$ von $M$ in allgemeiner Lage)
\begin{enumerate}
\item den \emph{Rang von $M$} als die Differenz der Dimensionen zwischen
einem $B$- und einem $N$-Orbit in $M$:
 \bdm
 \rk M\ :=\ \dim Bx - \dim Nx\, ;
 \edm
\item die \emph{Komplexit"at von $M$} als die Kodimension eines 
$B$-Orbits in $M$:
 \bdm
 c( M)\ :=\ \dim M - \dim Bx\ =\ \codim Bx\, .
 \edm
\end{enumerate}
Demnach ist die Kodimension eines allgemeinen $N$-Orbits genau $\rk M + c(X)$.
Nach dem Satz von Rosenlicht ist die Komplexit"at von $M$ ebenfalls gleich 
dem Transzendenzgrad von $\Rat{\C}{M}^B$ "uber $\C$; folglich sind 
die $G$-Variet"aten der Komplexit"at Null genau die sph"arischen 
$G$-Variet"aten, das sind diejenigen, f"ur die die Zerlegung
des Koordinatenrings vielfachheitenfrei ist. F"ur $c(M)=1$  
braucht der Invariantenring $\Aff{\C}{M}^G$ bereits nicht mehr trivial
zu sein (Beispiel: die $\SO(n,\C)$-Wirkung auf $\C^n$).
Weiterhin zeigte D.\ Panyushev in \cite[Cor.\ 2, p.253]{Panyushev90},
da"s $\rk M$ gleich der Dimension des von $\Spec$ erzeugten Teilraums in 
$P_+(G)\ox_{\Z}\Q$  und (mit Hilfe des Rosenlicht-Satzes) das 
Spektrum eine endlich erzeugte Halbgruppe ist.

Nimmt man $G$ als zusammenh"angend an, so ist der Fall $\rk M=0$ g"anzlich 
uninteressant, denn er beschreibt die
triviale Wirkung auf $M$. Wie Friedrich Knop bewiesen hat, mi"st
der Rang von $M$ -- siehe den nachstehenden Satz -- die  Gr"o"se des 
Zentrums von $\Diff[G]{M}$ (ebenso wie 
der Rang einer Lie-Algebra die Gr"o"se des Zentrums von 
$\env[\lalg]$ mi"st). Im glatten Fall entspricht die Komplexit"at in gewisser 
Weise der Anzahl positiver Wurzeln $|\Phi^+|=\dim \lalg[n]$, da 
f"ur den Transzendenzgrad von
$\Diff[G]{M}$ nach \cite[Satz 7.1]{Knop90} gilt
 \bdm
 \trdeg \Aff{\C}{T^*M}^G\ =\ 2c(M)+\rk M\, .
 \edm
%
\end{dfn}
\begin{thm}[Harish-Chandra-Homomorphismus {\cite[Thm.]{Knop94}}]
\label{harish-chandra}
Die zusammenh"angende reduktive komplexe algebraische Gruppe $G$ wirke
auf der  glatten affinen Variet"at $M$.  Dann existiert ein
Unterraum $\lalg[h]^*_M$ von $\lalg[h]^*$ der Dimension $\rk M$, eine
Untergruppe $W_M$ der Weylgruppe $W$ von $G$ (genannt die 
"`kleine Weylgruppe"' von $M$)
sowie ein Isomorphismus $\eta$ derart, da"s das folgende Diagramm
kommutiert:
 \bdm\begin{CD}
 \cent[\lalg[g]]@>>>\cent[M]\\
 @VV{\cong}V   @V{\eta}V{\cong}V\\
 \Aff{\C}{\lalg[h]^*}^W@>>\mathrm{res}>\Aff{\C}{\vrho+\lalg[h]^*_M}^{W_M}
 \end{CD}\edm
Insbesondere ist $\cent[M]$ als $\cent[\lalg[g]]$-Modul endlich erzeugt.
Zudem wirkt $W_M$ auf $\Aff{\C}{\vrho+\lalg[h]^*_M}$ als Spiegelungsgruppe,
d.h.\   $\cent[M]$  ist ein Polynomring in $\rk M$ Variablen mit homogenen
Erzeugenden.
\qed
\end{thm}\noindent
Die letzte Behauptung ist eine Folge aus dem Charakterisierungssatz f"ur
Pseudospiegelungsgruppen von Shepard-Todd und Chevalley 
(vgl.~\cite{Shepard&T54} und \cite{Chevalley55}).
%
Unter den Voraussetzungen von Korollar \ref{isotyp-zerl-aff} bezeichne
$U_{\lambda}=\Aff{\C}{M}^N (\lambda)$ f"ur $\lambda\in\Spec$
den Vielfachheitenraum zu $V^{\lambda}$, aufgefa"st als irreduzibler
$\Diff[G]{M}$-Modul. Weil nun $\cent[M]$ zu einem Polynomring in $n=\rk M$
Variablen isomorph ist, ist der zentrale Charakter $\chi_{\lambda}$ 
demnach die Auswertung eines Polynoms $\eta(Z)$ auf einem Element 
$\lambda'+\vrho\in\vrho+\lalg[h]^*_M$, genauer
 \bdm
 \chi_{\lambda}(Z)\ =\ \eta(Z)(\lambda'+\vrho),\quad Z\in\cent[M]\, .
 \edm
Die Frage, ob $U_{\lambda}$ durch $\lambda'$ eindeutig charakterisiert 
ist, haben wir bereits in Satz~\ref{central-char-uniq} positiv beantworten
k"onnen. Gemeinsam mit Satz~\ref{harish-chandra} bedeutet er damit, da"s 
statt der $\rk G$ Casimir-Operatoren bereits $\rk M$ zentrale 
Differentialoperatoren ausreichen, um das Spektrum von $M$ zu trennen.
Weiterhin erlaubt uns die Beschreibung von $\cent[M]$ als 
Polynomring zun"achst einzusehen, da"s
die Vielfachheitenr"aume aus der Frobenius-Zerlegung nur in Ausnahmef"allen
treue $\Diff[G]{M}$-Darstellungen liefern.
Zu einem $\Diff[G]{M}$-Modul $U$ und einer Teilmenge $X\subset U$ sei der
\emph{Annulator} definiert als
 \bdm
 \Ann (X)\ :=\ \{ D\in \Diff[G]{M}\ |\ Dx\ =\ 0\  \qqs x\in X\}\, .
 \edm
Ist $U$ zudem irreduzibel und $u$ ein beliebiges Element von $U-\{0\}$,
dann folgt aus dem Homomorphie-Satz sofort, da"s  $U$ zu 
$\Diff[G]{M}/ \Ann (u)$ isomorph ist, da $U$ in dieser Situation von
$u$  erzeugt wird.
\begin{kor}[Treue von {$\Diff[G]{M}$-Moduln}]\label{treue-invop-moduln}
Die zusammenh"angende reduktive komplexe algebraische Gruppe $G$ wirke
auf der  glatten affinen Variet"at $M$. 

Dann ist der Vielfachheitenraum  $U_{\lambda}$ genau dann
treuer $\Diff[G]{M}$-Modul , wenn $\rk M=0$ ist, 
also $G$ auf $M$ nur trivial wirkt.
%
\end{kor}
\begin{proof}
Wenn $G$ auf $M$ trivial wirkt, so ist es klar, da"s $\Diff{M}$ auf
$\Aff{\C}{M}$ treu wirkt. Ist dagegen $\rk M =:n>0$, so bemerken wir,
da"s nat"urlich
 \begin{equation}\label{ann-diff-mod}
 \Diff[G]{M} \cdot \ker \chi_{\lambda}\ \subset\ \Ann (U_{\lambda})
 \end{equation}
gilt; damit reicht es zu zeigen, da"s $\ker \chi_{\lambda}$ nicht trivial ist.
Aber in diesem Fall ist $\cent[M]$
isomorph zu einem Polynomring in $n$ Variablen, der zentrale
Charakter $\chi_{\lambda}$ demnach ein (multiplikativer) Homomorphismus von 
$\Aff{\C}{X_1,\ldots,X_n}$ nach $\C$. Ein solcher ist bereits durch die Bilder
$c_i$ der Unbestimmten $X_i$ eindeutig festgelegt, und jedes Polynom, welches
einen Linearfaktor $X_i-c_i$ enth"alt, liegt im Kern von $\chi_{\lambda}$.
\end{proof}

\medskip\noindent Es scheint nicht bekannt zu sein, ob 
in Gleichung~\ref{ann-diff-mod} -- wie im Falle
der Verma-Moduln \cite[Thm.8.4.3]{Dixmier2} -- statt der Inklusion
die Gleichheit gilt.
%
\section{Die Frobenius-Zerlegung f"ur reelle Formen}
\noindent
In der Praxis tritt h"aufig der Fall auf, da"s die komplexe Situation
in Wahrheit durch Komplexifizierung einer reellen Gruppenwirkung entsteht.
Deswegen wollen wir hier die Korrespondenz zwischen reeller und 
komplexer Frobenius-Zerlegung herstellen.
\begin{dfn}
V"ollig kanonisch l"a"st sich zu jeder reellen affinen Variet"at
$M_0\subset \R^n$ deren Komplexifizierung $M_0(\C)=:M$ in $\C^n$ konstruieren
als
 \bdm
 M:= \{ x\in\C^n\ | \ f(x) = 0 \ \qqs f\in \I (M_0)\},
 \edm
wobei $\I (M_0)$ das Verschwindungsideal von $M_0$ bezeichne.
Eine \emph{reelle Form} der komplexen affinen Variet"at $M$  soll nun eine
reelle affine Variet"at $M_0$ sein, deren Komplexifizierung $M_0(\C)$ zu $M$
isomorph ist.
In diesem Fall liegt $M_0$ dicht in $M$; weil eine affine Variet"at genau dann 
irreduzibel ist, wenn ihr Abschlu"s es ist, folgt damit, da"s $M_0$ genau
dann irreduzibel ist, wenn $M$ es ist. (vgl.\  \cite[S.61-64]{Onishchik&V2}).
Es existiert ein Vektorraum-Isomorphismus zwischen den
komplexwertigen regul"aren Funktionen auf $M_0$ und den "ublichen regul"aren
Funktionen auf $M$
 \bdm
 \Aff{\R}{M_0}\ox_{\R}\C\ \cong\  \Aff{\C}{M};
 \edm
denn in die eine Richtung ist dies einfach die Einschr"ankung einer
Funktion von $M$ auf $M_0$; in der anderen Richtung 
impliziert die Dichtheit von $M_0$ in $M$, da"s jede regul"are
Funktion auf $M_0$ eine eindeutige regul"are Fortsetzung auf ganz $M$ hat.
Deswegen schreibt man oft $\Aff{\C}{M}$, wenn man $\Aff{\R}{M_0}\ox\C$ meint.
\end{dfn}
\noindent
Damit gilt insbesondere f"ur eine glatte reelle affine Variet"at $M_0$ mit 
Komplexifizierung $M$,
da"s $\Diff{M_0}\ox\C \cong \Diff{M}$ (zun"achst im Vektorraumsinne) gilt.
In der Tat, ist $D_0\in \Diff{M_0}$ ein Differentialoperator auf $M_0$,
so l"a"st sich eine Fortsetzung $D$ von $D_0$ auf ganz $M$ eindeutig definieren
durch die Forderung, da"s das Diagramm
 \bdm\begin{CD}
 \Aff{\R}{M_0}_{\C}@>{\cong}>>\Aff{\C}{M}\\
 @V{D_0}VV   @VV{D}V\\
 \Aff{\R}{M_0}_{\C}@>{\cong}>>\Aff{\C}{M}\\
 \end{CD}\edm
kommutieren m"oge, und es macht gleichzeitig klar, was unter der Einschr"ankung
eines Operators $D\in\Diff{M}$ auf $M_0$ zu verstehen sei.
Klarerweise ist die Fortsetzung / Einschr"ankung
von Differentialoperatoren ein Algebrenisomorphismus. 
Weil wir uns im Reellen meist f"ur komplexwertige regul"are Funktionen 
interessieren,
macht es Sinn, von nun an die entsprechenden Differentialoperatoren
$\Diff{M_0}_{\C}$ zu schreiben. Es folgt sofort, da"s ein komplexer 
$\Diff{M_0}_{\C}$-Modul genau dann
irreduzibel ist, wenn er als $\Diff{M}$-Modul irreduzibel ist.
\begin{dfn}
Eine reelle algebraische Gruppe $G_0$ soll \emph{reelle Form} der komplexen
algebraischen Gruppe $G$ genannt werden, wenn sie eine reelle Form im eben 
genannten Sinne affiner Variet"aten ist und die identische Einbettung
$G_0\inj G$ sich zu einen Gruppenisomorphismus $G_0(\C)\cong G$ fortsetzen
l"a"st. Zum Beispiel ist, "ahnlich wie bei Lie-Gruppen,  die Fixpunktmenge 
eines beliebigen involutiven 
antiholomorphen Automorphismus einer komplexen zusammenh"angenden Gruppe
 $G$ eine reelle Form (vgl.\  \cite[S.100-102]{Onishchik&V2}). 
Man beweist dann leicht folgendes Lemma:
\begin{lem}[Komplexe Darstellungen reeller Formen reduktiver Gruppen]
\label{complexreps-realforms}
Sei $G_0\subset G$ eine irreduzible reduktive reelle algebraische Gruppe und
zudem eine reelle
Form der (ebenfalls irreduziblen und reduktiven) komplexen algebraischen
Gruppe $G$.
Dann gibt es eine Bijektion zwischen den komplexen irreduziblen regul"aren
Darstellungen von $G_0$ und denen von $G$,
welche durch Restriktion bzw.\ eindeutige Fortsetzung gegeben ist. \qed
\end{lem}
\begin{proof}
Sei $(\vrho_0,V)$ eine (regul"are endlich-dimensionale) irreduzible 
Darstellung von $G_0$ auf dem komplexen Vektorraum $V$, also ein Morphismus
$G_0\ra \GL(V)$. W"ahle nun eine Basis $e_1,\ldots,e_n$ von $V$. 
Dann lassen sich die Funktionen 
 \bdm
 m_{ij}:\quad G_0\lra G,\quad \vrho_0(g)\ =\ \sum m_{ij}(g)e_i
 \edm
als Elemente von $\Aff{\R}{G_0}_{\C}=\Aff{\C}{G}\restr_{G_0}$ eindeutig
auf ganz $G$ fortsetzen, weil $G_0$ in $G$ dicht liegt. Diese definieren
somit eine Darstellung von $G$ auf $V$, die klarerweise irreduzibel ist. 
Umgekehrt ist die Restriktion einer irreduziblen Darstellung 
$(\vrho,V)$ von $G$ auf $G_0$ nat"urlich wieder eine $G_0$-Darstellung; einzig
nicht trivial an diesem Lemma ist folglich der Nachweis, da"s diese
Einschr"ankung irreduzibel ist, wof"ur wir auf eine Idee
von \v{Z}elobenko (vgl.\ \cite[Ch.VI, \S 42, Thm.2]{Zelobenko1}) zur"uckgreifen
wollen. 

Sei $V_0\neq \{0\}$ ein $G_0$-invarianter Unterraum von $V$. Angenommen,
$V_0$ ist nicht $G$-invariant: w"ahle  $g\in G$ und $v_0\in V_0$ mit
$\vrho(g)v_0\sla{\in} V_0$. Sodann existiert ein lineares Funktional
$f:V\ra \C$ derart, da"s $f\restr_{V_0}=0$ und $f(\vrho(g)v_0)\neq 0$ gilt. 
Ein Dichtheitsargument zeigt, da"s es ein solches nicht geben kann:
denn ist $f:V\ra\C$ ein beliebiges lineares Funktional, welches auf $V_0$
verschwindet, dann betrachte man zu jedem $v_0\in V_0$ die -- offensichtlich
regul"are -- Abbildung $f_{v_0}:G_0\ra\C$, $g_0\mapsto f(\vrho(g_0) v_0)=0$.
Weil $G_0$ in $G$ dicht ist, hat diese Abbildung nur die Nullabbildung als
regul"are Fortsetzung auf ganz $G$; das hei"st aber nichts anderes, als
da"s $f(\vrho(g) v_0)=0$ ist f"ur alle $g\in G$ und alle $v_0\in V_0$, im
Widerspruch zum weiter oben
konstruierten Funktional. Also ist $V_0$ doch $G$-invariant, damit 
$V_0=V$ und $(\vrho\restr_{G_0},V)$ wie gew"unscht $G_0$-irreduzibel.
\end{proof}
\end{dfn}
\begin{dfn}
Sei $(G,M)$ ein Paar bestehend aus einer komplexen algebraischen reduktiven
zusammenh"angenden Gruppe $G$, einer ebenfalls komplexen affinen glatten 
irreduziblen Variet"at $M$, sowie 
einer regul"aren Wirkung von $G$ auf $M$. Wir nennen ein Paar $(G_0,M_0)$,
bestehend aus einer reellen algebraischen Gruppe $G_0$ und einer affinen 
glatten reellen Variet"at $M_0$, eine \emph{reelle Form}  von $(G,M)$, wenn
$G$ bzw.\ $M$ die Komplexifizierung von $G_0$ bzw.\ $M_0$  und die
$G$-Wirkung auf $M$ ebenfalls die Komplexifizierung der $G_0$-Wirkung auf 
$M_0$ ist.
\end{dfn}
\noindent
Sei nun $(G_0,M_0)$ eine reelle Form von $(G,M)$, und weiterhin
 \bdm
 r:\  \Aff{\C}{M}\ \stackrel{\cong}{\lra}\ \Aff{\R}{M_0}_{\C}
 \edm
der Algebren-Isomorphismus zwischen den ("ublichen) regul"aren Funktionen auf
$M$ und den komplexwertigen regul"aren Funktionen auf $M_0$, der durch
Restriktion bzw. eindeutige Fortsetzung von Funktionen gegeben ist.
In dieser Situation liefert Korollar~\ref{isotyp-zerl-aff} die 
Frobenius-Zerlegung von $\Aff{\C}{M}$.
Vergessen wir einen Augenblick die unterliegenden Variet"aten und 
betrachten nur die
Wirkung $\vrho$ von $G$ auf dem komplexen Vektorraum $\Aff{\C}{M}$.
Nach Lemma~\ref{complexreps-realforms} ist dann die isotypische
Zerlegung von $\Aff{\C}{M}$ unter $G$ identisch mit der isotypischen
Zerlegung von $\Aff{\C}{M}$ unter der Wirkung von $G_0$, definiert als die
Einschr"ankung der $G$-Wirkung hierauf, geschrieben $\vrho_0$. Die
Einschr"ankung auf $G_0$ der komplexen $G$-Darstellung $V^{\lambda}$
schreiben wir  $V^{\lambda}_0$. Nun verwenden wir, da"s
$\Aff{\C}{M}$ mit der Abbildung $r$ isomorph ist zu $\Aff{\R}{M_0}_{\C}$, die 
Einschr"ankung der $G_0$-Wirkung nat"urlich mit der "Ubertragung der
$G_0$-Wirkung von $\Aff{\R}{M_0}$ auf $\Aff{\R}{M_0}_{\C}$ "ubereinstimmt
und folglich Korollar \ref{isotyp-zerl-aff} ebenfalls die isotypische
Zerlegung von $\Aff{\R}{M_0}_{\C}$ unter $G_0$ liefert. Insgesamt
haben wir also als $G_0$-Modul  -- zun"achst vom abstrakten Standpunkt
der Morphismen -- die Zerlegung
 \bdm
 \Aff{\R}{M_0}_{\C}\ \cong\  \bigop_{\lambda\in P_+(G)} 
 \Hom_{G_0}(V^{\lambda}_0,\Aff{\R}{M_0}_{\C} )\ox V^{\lambda}_0\, ,
 \edm
welcher auf der Seite der Funktionen unter Verwendung des Isomorphismus $r$
die Summe
 \bdm
 \Aff{\R}{M_0}_{\C}\ \cong\  \bigop_{\lambda\in P_+(G)} 
 r\left(\Aff{\C}{M}^N(\lambda)\right)\ox V^{\lambda}_0
 \edm
entspricht. Im allgemeinen ist dies die einzige Beschreibung der reellen
Frobenius-Zerlegung in der Sprache des Transformationsverhaltens gewisser
Funktionen. Eine intrinsische Charakterisierung der Vielfachheitenr"aume
alleine aus der Struktur der reellen Gruppe $G_0$ ist nur in einem
Ausnahmefall m"oglich, denn wir nun skizzieren wollen.

\medskip\noindent
Sei die reelle Form $G_0$ definiert als die Fixpunktmenge des involutiven
Antiautomorphismus $\sigma$, und $(\tilde{G}_0,\tau)$ eine damit vertr"agliche 
reelle kompakte Form von $G$, d.h.\ es gelte $\sigma\tau=\tau\sigma$.
Dann ist $\theta:=\sigma\tau$ die Cartan-Involution von $G_0$,
die Lie-Algebra $\lalg_0$
zerlegt sich in die Eigenr"aume zum Eigenwert $+1$ und $-1$ 
 \bdm
 \lalg_0\ =\ \lalg_0(1) \op\lalg_0(-1)\ =:\ \lalg[k]\op\lalg[p],  
 \edm
und es existieren immer Cartan-Unteralgebren von $\lalg_0$, welche unter
$\theta$ invariant ist. Sei fortan $\lalg[h]_0$ eine solche, und $H_0$ eine 
Untergruppe von $G_0$ mit Lie-Algebra $\lalg[h]_0$.
Es ist $H_0$ eine reelle Form von $H$, der Untergruppe von 
$G$ mit Lie-Algebra $\lalg[h]:=\lalg[h]_0(\C)$; deswegen liegt $H_0$ in $H$ 
dicht. Ist nun  $B=HN$ die
Erg"anzung von $H$ zu einer Boreluntergruppe von $G$, die durch
die Wurzelzerlegung bez"uglich $H$ definiert ist, so braucht $N$
nicht unter $\sigma$ invariant zu sein, deswegen bilden die  Fixpunkte von 
$N$ unter $\sigma$ in $G_0$ im allgemeinen keine reelle Form von $N$.
Wei"s man aber, da"s $\lalg[h]_0$ vollst"andig in $\lalg[p]$ enthalten ist, 
so ist $\lalg[h]_0$ gleichzeitig eine maximal abelsche Unteralgebra von
$\lalg[g]_0$ und f"ur alle  $h\in\lalg[h]_0$ ist der Operator $\ad h$ 
reell diagonalisierbar. Daraus folgt, da"s die 
Wurzelr"aume $\lalg^{\alpha}$ der komplexen Lie-Algebra genau die
Komplexifizierungen der  -- v"ollig
analog definierten -- reellen Wurzelr"aume ($\alpha\in\lalg[h]_0^*$)
 \bdm
 \lalg[g]^{\alpha}_0\ :=\ \{ x\in\lalg_0 : \ad\, h (x)= \alpha(h)x\}
 \edm
sind; bezeichne $\Delta_+$ die positiven Wurzeln (von $\lalg$ oder $\lalg_0$,
dies ist nun egal), so ist insbesondere  
 \bdm
 \lalg[n]_0\ :=\ \sum_{\alpha\in\Delta_+} \lalg[g]^{\alpha}_0\
 \edm
eine reelle Form von $\lalg[n]$, der Lie-Algebra von $N$. Die Fixpunkte
$N_0$ von $N$ unter $\sigma$ liegen also dicht in $N$. 
Reelle Lie-Algebren, die eine vollst"andig in $\lalg[p]$ enthaltene
Cartan-Unteralgebra besitzen, hei"sen \emph{zerfallend} ("`split"', 
"`d\'eploy\'ee"'). F"ur die klassischen einfachen Lie-Algebren sind dies
bekanntlich $\lalg[sl](n,\R)$ f"ur $n\geq 2$, $\lalg[so](k,k+1)$ f"ur 
$k\geq 1$,
$\lalg[so](k,k)$ f"ur $k\geq 3$ und $\lalg[sp](n,\R)$ f"ur $n\geq 2$.

F"ur die urspr"unglich
gestellte Aufgabe, die reelle Frobenius-Reziprozit"at alleine anhand der
reellen Gruppe $G_0$ zu charakterisieren, sind wir nun in einer
"au"serst komfortablen Situation; denn jetzt hat jede Funktion in
$\Aff{\R}{M_0}_{\C}$, die sich  unter $B_0:=H_0N_0$ gem"a"s dem Gewicht
$\lambda$ transformiert, eine eindeutige Fortsetzung zu einer Funktion aus
$\Aff{\C}{M}$ mit gleichen Transformationseigenschaften bzgl.\ $B$,
und demnach ist
 \bdm
 \Aff{\R}{M_0}_{\C}\ \cong\  \bigop_{\lambda\in P_+(G)} 
 \Aff{\R}{M_0}^{N_0}_{\C}(\lambda)\ox V^{\lambda}_0\, .
 \edm
Dieses Resultat kann man sogar noch dahingehend verbessern, da"s
man sich auf reellwertige Funktionen beschr"anken kann: denn leicht "uberlegt
man sich, da"s f"ur eine zerfallende Lie-Algebra $\lalg[g]_0$
die komplexen Darstellungen von $\lalg[g]_0(\C)=:\lalg$ 
genau die Komplexifizierungen der reellen Darstellungen von $\lalg_0$
sind. Wir fassen die Ergebnisse dieses Abschnitts zusammen:
\begin{kor}[Reelle Frobenius-Zerlegung]
Sei $G_0$ eine reelle reduktive zusammenh"angende algebraische Gruppe,
die reelle Form der komplexen reduktiven algebraischen Gruppe $G$ ist, und 
$M_0\subset \R^n$ eine reelle, glatte, irreduzible Variet"at
mit Komplexifizierung $M$, auf der 
$G_0$ regul"ar wirkt. Sei weiterhin 
$r:\ \Aff{\C}{M}\ \stackrel{\cong}{\lra}\ \Aff{\R}{M_0}_{\C}$ der
zugeh"orige Algebren-Isomorphismus. Dann gilt:
\bdm
 \Aff{\R}{M_0}_{\C}\, \cong\,  \bigop_{\lambda\in P_+(G)} 
 \Hom_{G_0}(V^{\lambda}_0,\Aff{\R}{M_0}_{\C} )\ox V^{\lambda}_0\, 
 \cong\,  \bigop_{\lambda\in P_+(G)} 
 r\left(\Aff{\C}{M}^N(\lambda)\right)\ox V^{\lambda}_0
\edm
als $G_0$-Moduln, und, wenn $G_0$ zus"atzlich entfaltet ist,
\bdm
\Aff{\R}{M_0}\, \cong\,  \bigop_{\lambda\in P_+(G)} 
 \Aff{\R}{M_0}^{N_0}(\lambda)\ox V^{\lambda}_0\, .
\edm
Zudem sind die Vielfachheitenr"aume der $V^{\lambda}_0$ als
$\Diff[G]{M}_{\C}$-Moduln irreduzibel, paarweise in"aquivalent und
haben zentralen Charakter.
\qed
%
\end{kor}
\section{Einige singul"are Beispiele}\label{sing-exa}\noindent
Die Glattheit der zugrunde liegenden Variet"at $M$ wurde im Beweis der 
allgemeinen Frobenius-Rezoprozit"at an genau zwei Stellen verwendet. Zum
einen garantiert sie wegen Lemma \ref{diff-simple-alg} die
G"ultigkeit des Ergebnisses im Falle der trivialen Gruppenwirkung.
Eben diese Eigenschaft braucht man zudem, um den
Dichtheitssatz von Jacobson
zur Darstellbarkeit von Endomorphismen endlich-dimensionaler Teilr"aume
als Differentialoperatoren sowie das Lemma von Dixmier anwenden
zu k"onnen (dieses wird jedoch erst beim Studium der zentralen 
Charaktere ben"utzt).
Zum anderen erlaubt erst die Identifizierung von $\Diff{M}$ mit 
$\Aff{\C}{T^*M}$ "uber die Symbolabbildung den Schlu"s, da"s die
Wirkung der Gruppe $G$ auf den algebraischen Differentialoperatoren 
lokal endlich ist; hat man nur einen Unterring von $\Aff{\C}{T^*M}$, so
w"are denkbar, da"s beim "Ubergang zu $G$-invarianten Teilr"aumen, die einen
vorgegebenen kleineren Teilraum enthalten,
Funktionen entstehen, die nicht Symbol eines Differentialoperators sind. 

Schlu"sendlich ist es zwar nicht f"ur den Beweis, wohl aber konzeptionell
und f"ur das
konkrete Rechnen st"orend, da"s man eine Wahl treffen mu"s zwischen der 
Algebra aller Differentialoperatoren $\Diff{M}$ und denen, die nur von den
Funktionen und den Derivationen von $M$ erzeugt werden, $\Diff{M}^*$
in der Notation von Abschnitt \ref{section-alg-diff-ops}. Die erste 
Algebra l"a"st sich zwar definieren,
aber nicht bestimmen, kann beliebig unanst"andige Eigenschaften
haben und liefert keine Differentialoperatoren im Sinne der
Analysis; die zweite Algebra ist handlicher und  analytisch sinnvoll, aber 
meist zu "`klein"'.  

Gegenstand dieses Abschnitts ist es deshalb, anhand einiger Beispiele
diejenigen Effekte zu illustrieren, die bei reduktiven Gruppenwirkungen auf
singul"aren affinen Variet"aten auftreten k"onnen. Sie zeigen, da"s allgemeine
S"atze hier nicht zu erwarten sind. Die Algebra der
Wahl ist dabei im Grunde immer $\Diff{M}^*$; "uber $\Diff{M}$ kann man nur
verzeinzelt Aussagen machen, wenn man irgendwie durch eine g"unstige Wahl der Koordinaten Elemente hieraus "`raten"' kann.

\begin{exa}[Spitzer Doppelkegel]
Wir betrachten den spitzen Doppelkegel im $\R^3$
 \bdm
 X_0\ :=\ \{(x,y,z)\in\R^3\ |\ x^2+y^2=z^2\}\, . 
 \edm
Er ist eine reelle Form des komplexen Kegels
 \bdm
 X \ :=\ \{(x,y,z)\in\C^3\ |\ x^2+y^2=z^2\}\, , 
 \edm
und der singul"are Ort besteht in beiden F"allen genau aus dem Nullpunkt.
Folglich k"onnen wir den Koordinatenring $\Aff{\C}{X}$ von $X$ mit den
komplexwertigen regul"aren Funktionen $\Aff{\R}{X_0}_{\C}$ auf $X_0$ 
identifizieren. Parametrisiert man $X_0$ mittels
 \bdm
 x\ =\ u\cos \theta,\quad y\ =\ u\sin \theta,\quad z\ =\ u,
 \edm
so wird nach einer komplexen Drehung der Ring $A:=\Aff{\R}{X_0}_{\C}$
"uber $\C$ erzeugt von $u$, $u e^{i\theta}$ und $u e^{-i\theta}$
 \bdm
 A\ :=\ \Aff{\R}{X_0}_{\C}\ \cong\ \Aff{\C}{u,u e^{i\theta},u e^{-i\theta}}\, .
 \edm
Bevor wir eine bestimmte Gruppenwirkung auf dem Kegel $X_0$ w"ahlen, bestimmen
wir alle seine Derivationen sowie die Wirkung der von ihnen und
$\Aff{\R}{X_0}_{\C}$ erzeugten Algebra $\Diff{X_0}^*$. Zun"achst ist klar,
da"s die Koordinatenlinien die Derivationen
 \bdm
 u\del_u\ =\ x\del_x+y\del_y+z\del_z\ \ \text{ und }\ \ 
 \del_{\theta}\ =\ x\del_y-y\del_x
 \edm
erzeugen. Die Algebra $\Der X_0$ enth"alt aber noch weit mehr Elemente:
\begin{lem}[Derivationen des Kegels]\label{deriv-kegel}
Durch folgende Kriterien ist $\Der X_0$ vollst"andig bestimmt:
\begin{enumerate}
\item Ein Element der Gestalt $q\del_u$ liegt genau dann in $\Der X_0$,
wenn $\,q\in u\cdot A$;
\item Ein Element der Gestalt $q\del_{\theta}$ liegt genau dann in $\Der X_0$,
wenn $\,q\in A$;
\item Ein Element der Gestalt $q_1\del_u + q_2\del_{\theta}$ liegt genau dann 
in $\Der X_0$, wenn es Summe zweier Derivationen der gerade genannten
Gestalt oder in der $A$-linearen H"ulle der beiden Derivationen
 \bdm
 R\ :=\ u e^{i\theta}\del_u + i e^{i\theta}\del_{\theta},\quad
 S\ :=\ u e^{-i\theta}\del_u - i e^{-i\theta}\del_{\theta}
 \edm
enthalten ist.
\end{enumerate}
\end{lem}
\begin{proof} Zum Beweis verwenden wir Das Kriterium aus 
Lemma \ref{diff-ops-krit}.
Wir setzen eine Derivation von $X_0$ an mit $D=q\del_u$, 
$q\in\Rat{\R}{X_0}_{\C}$. Wegen der Derivationseigenschaft gen"ugt es,
wenn die Bilder der Koordinatenfunktionen unter $D$ in $A$ liegen, denn dann 
folgt sofort $D(A)\subset A$. Es m"ussen also
 \bdm
 D(u)\ =\ q,\quad D(u e^{i\theta})\ =\ q e^{i\theta},\quad 
 D(u e^{-i\theta})\ =\ q e^{-i\theta}
 \edm
Elemente von $A$ sein, d.h.\ $q$ liegt in $A$ und, verm"oge Multiplikation
mit $u e^{-i\theta}$ der zweiten Bedingung, $q u$ liegt ebenfalls in $A$.
Schreibt man dann $q u $ als Summe
 \bdm
 qu\ =\ \sum \alpha\, u^n (u e^{i\theta})^m(u e^{-i\theta})^l, 
 \edm
so implizieren die erste und die zweite Bedingung, da"s in allen
auftretenden Summanden $m$ und $l$ gr"o"ser Eins sein m"ussen, d.h.
 \bdm
 qu\ =\ u^2 \sum \alpha u^n (u e^{i\theta})^{m-1}(u e^{-i\theta})^{l-1}\ =:\
 u^2\cdot a
 \edm
mit einem $a\in A$; es folgt $q=ua$ und damit hat $q$ die behauptete
Gestalt. Analog beweist man den zweiten Teil des Lemmas. F"ur den
letzten Teil machen wir den Ansatz $D=q_1\del_u + q_2\del_{\theta}$,
was impliziert, da"s 
 \bdm
 D(u)\ =\ q_1,\quad D(u e^{i\theta})\ =\ e^{i\theta}(q_1+iuq_2),\quad 
 D(u e^{-i\theta})\ =\ e^{-i\theta}(q_1-iuq_2)
 \edm
regul"are Funktionen sein sollen. Nun betrachtet man die $u$-Graduierung
von $A$: wenn $\deg q_1\neq \deg q_2+1$ ist, so haben die beiden Summanden
in der ersten und zweiten Bedingung unterschiedlichen $u$-Grad, m"ussen
also jeder f"ur sich in $A$ liegen, was einer Derivation der soeben
beschriebenen Gestalt entspricht. Sei deswegen fortan 
$\deg q_1= \deg q_2+1$. Insbesondere ist dann $\deg q_1\geq 1$, ich kann
also $q_1=uq_1'$ setzen. Aus der zweiten Bedingung folgt mittels
Multiplikation mit $u e^{i\theta}$, da"s $u^2(q_1'+iq_2)\in A$ sein mu"s;
aus der ersten Bedingung $q_1=q_1'u\in A$ folgt aber $q_1'u^2\in A$,
also durch Subtraktion $u^2q_2\in A$. Man "uberlegt sich leicht, da"s
hieraus folgt, da"s $q_2$ die Gestalt
 \bdm
 q_2\ =\ a_2 + b_2 e^{i\theta} + b_2' e^{-i\theta}+ c_2 e^{2i\theta} + 
 c_2' e^{-2i\theta}\, , 
 \edm
haben mu"s, wobei $a_2$ in $A$ ist und man die anderen Koeffizienten so w"ahlen
kann, da"s sie zwar in $A$, aber $b_2$, $b_2'$ nicht in $u\cdot A$ und 
$c_2$, $c_2'$ nicht in 
$u^2\cdot A$ liegen. In "ahnlicher Weise kann man $q_1'$ darstellen.
Weil aber hier bereits $uq_1'$ in $A$ sein mu"s, verschwinden die
beiden letzten Terme:
 \bdm
 q_1\ =\ a_1 + b_1 e^{i\theta} + b_1' e^{-i\theta}\, .
 \edm
Schreibt man nun das in der zweiten Bedingung vorkommende Element 
$u e^{i\theta}(q_1'+iq_2)$ in dieser Darstellung 
aus und subtrahiert alle Elemente, die trivialerweise in $A$ liegen, so bleibt,
da"s $ue^{2i\theta}(b_1+ib_2)+ic_2e^{i3\theta}u$ in $A$ liegen mu"s, was
nach Wahl der Koeffizienten nicht geht; deswegen ist $c_2=0$ und
$b_1+ib_2=0$. Ebenso zeigt man $c_2'=0$ und $b_1'-ib_2'=0$. Die
"`Minimall"osungen"' $b_1=1$, $b_2=i$ und $b_1'=1$, $b_2'=-i$ entsprechen
den Elementen $R$ und $S$. F"ur diese berechnen wir f"ur sp"ater noch
die Produkte
 \begin{equation}\label{rel-RS}
 RS\ =\  u^2\del^2_u + 2 u\del_u+\del^2_{\theta}-i\del_{\theta},\quad
 SR\ =\  u^2\del^2_u + 2 u\del_u+\del^2_{\theta}+i\del_{\theta},\quad
[R,S]\ =\ -2i\del_\theta\, .
 \end{equation}
\end{proof}
%
Zur Illustration, wie die Derivationen des Kegels auf seinen regul"aren
Funktionen operieren, dient Abbildung \ref{kegel-koordring}. Darin
sind alle regul"aren Funktionen in einem rechteckigen  
Zeilen- und Spaltenmuster  angeordnet; die positive $x$-Achse beschreibt
die wachsenden $u$-Potenzen, die $y$-Achse die ganzzahligen Vielfachen
von $\theta$ im Argument der Exponentialfunktion. Wie man sieht, wird von 
allen m"oglichen Besetzungsstellen nur ein (gedrehter) Quadrant tats"achlich
belegt: dies liegt daran, da"s z.\ B.\ zwar $ue^{i\theta}$ in $A$ liegt, nicht
aber $ue^{2i\theta}$. Die von $A$ und $\Der X_0$ erzeugte Algebra 
$\Diff{X_0}^*$ enth"alt zun"achst die Multiplikation mit allen 
Koordinatenfunktionen; beginnt man bei der konstanten Funktion Eins
im linken Bildteil, so sind dies genau die drei nach rechts zeigenden
Pfeile (zwei diagonal, einer horizontal). von jeder anderen Funktion gehen
nat"urlich immer diese drei Pfeile aus; der "Ubersichtlichkeit halber
haben wir sie nur am oberen und unteren Bildrand eingezeichnet, und ansonsten
durch ein Sternchen angedeutet (jedes solche steht also f"ur ein nach
rechts zeigendes Pfeilkreuz). Interessant ist nun die Wirkung der
Derivationen: $u\del_u$ und $\del_{\theta}$ sind Multiplikationen mit einer
Konstanten und deswegen nicht eingezeichnet. Die einzig nichttrivialen
Derivationen sind die Elemente $R$ und $S$. Sie entsprechen
(bis auf Vielfache, die hier nicht interessieren) einem "`Auf-"' und 
einem "`Absteiger"' bei fixiertem $u$-Grad.
 
\begin{figure}
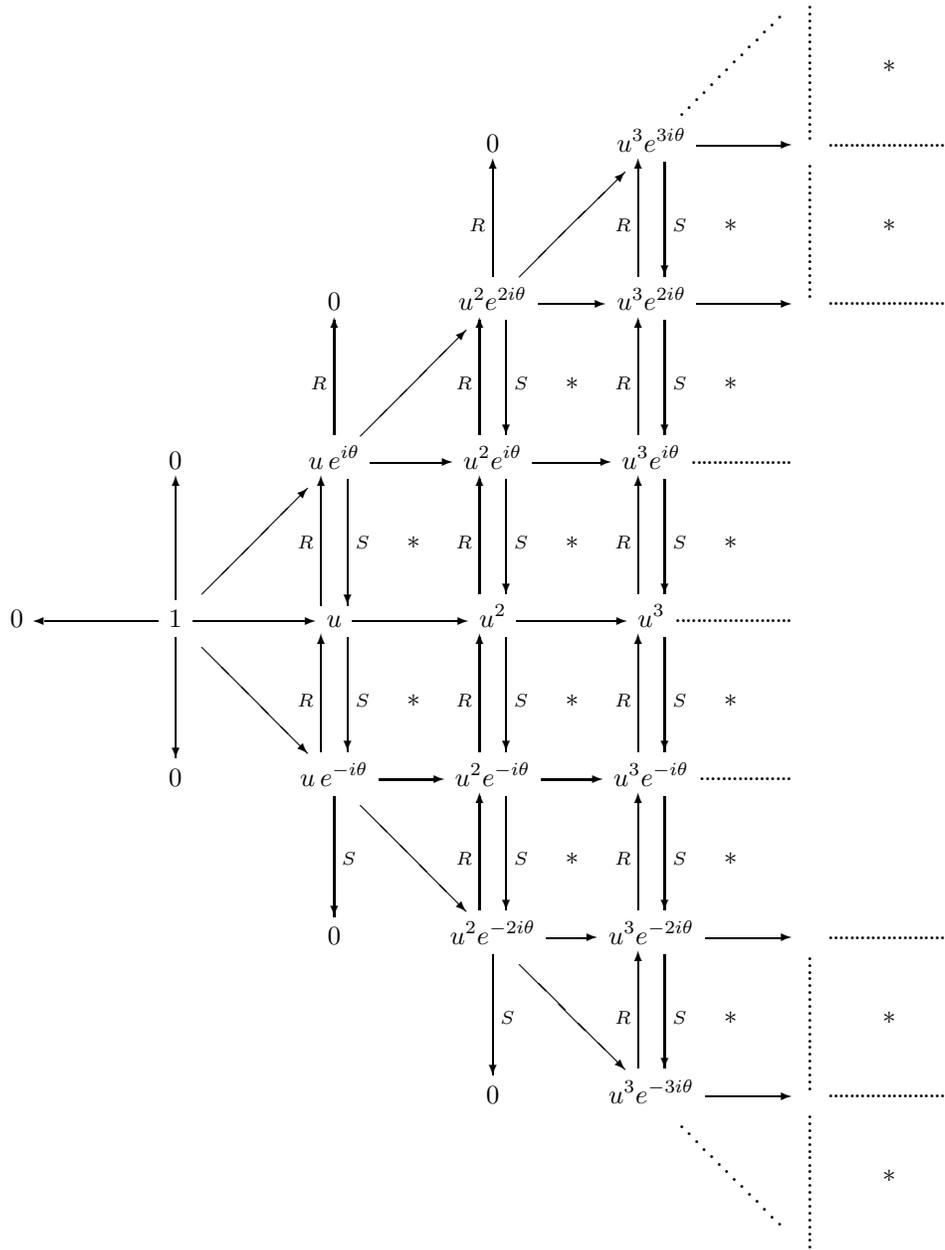

\begin{diagram}[labelstyle=\scriptstyle]
 & & & & & & & & & & \hphantom{x} & & \\
 & & & & & & & & & \ruDots & \uDots & * & \\
 & & & & & & 0 & & u^3 e^{3i\theta} & \rTo &\hphantom{x} &\rDots & 
    \hphantom{x}\\
 & & & & & & \uTo^R & \ruTo &\uTo^R \dTo_S & * &\uDots & * & \\
 & & & & 0 & & u^2 e^{2i\theta} &\rTo & u^3 e^{2i\theta}& \rTo & 
    \hphantom{x}& \rDots & \hphantom{x} \\
 & & & & \uTo^R & \ruTo& \uTo^R\dTo_S & * &\uTo^R\dTo_S & *  & &  \\
 & & 0 & & u\, e^{i\theta} &\rTo & u^2 e^{i\theta}  & \rTo & 
   u^3 e^{i\theta}&  \rDots & \hphantom{x}& & \\
 & & \uTo & \ruTo & \uTo^R\dTo_S & * & \uTo^R\dTo_S & * &\uTo^R\dTo_S &
    * & & \\
0 & \lTo & 1 & \rTo & u  &\rTo & u^2   & \rTo & u^3 & \rDots & 
    \hphantom{x}& &\\
 & & \dTo & \rdTo & \uTo^R\dTo_S & * & \uTo^R\dTo_S &* &\uTo^R\dTo_S &
    * & & & \\
 & & 0 & & u\, e^{-i\theta} &\rTo & u^2 e^{-i\theta}  & \rTo & 
   u^3 e^{-i\theta}& \rDots &\hphantom{x} & & \\
 & & & & \dTo_S & \rdTo & \uTo^R\dTo_S & * &\uTo^R\dTo_S & * & & & \\
 & & & & 0 & & u^2 e^{-2i\theta} &\rTo & u^3 e^{-2i\theta}& \rTo & 
    \hphantom{x}& \rDots & \hphantom{x} \\
 & & & & & & \dTo_S & \rdTo &\uTo^R\dTo_S & * & \uDots & * & \\
 & & & & & & 0 & & u^3 e^{-3i\theta} & \rTo & \hphantom{x} & \rDots &
     \hphantom{x}\\
 & & & & & & & & & \rdDots & \dDots& * & \\
 & & & & & & & & & & \hphantom{x}& & \\
\end{diagram}
\caption{Koordinatenring des Doppelkegels $z^2=x^2+y^2$}
\label{kegel-koordring}
\end{figure}
Als m"ogliche Gruppenwirkungen auf $X_0$ bieten sich nun drei Kandidaten
an. Die kanonische Wirkung der $\SO{}^0(2,1)$ ist relativ uninteressant,
weil vielfachheitenfrei; ihre eindimensionalen Vielfachheitenr"aume
sind wenig geeignet, um irgendwelche Effekte zu demonstrieren
(es sei bemerkt, da"s die Derivationen $R$ und $S$ eigentlich
von dieser Gruppenwirkung stammen). Stattdessen betrachten wir die
multiplikative Wirkung von $\R^*$ auf $X_0$ (von $\C^*$ auf $X$),
sowie die Drehwirkung der $\SO(2,\R)$ in jeder $x$-$y$-Ebene bei 
festgehaltener $z$-Koordinate. Diese Gruppenwirkungen haben genau
die "`Spalten"' bzw.\ die "`Zeilen"' von $A$ als Vielfachheitenr"aume.
Es gilt:
%
\begin{prop}
\begin{enumerate}
\item Die Vielfachheitenr"aume der $\R^*$-Wirkung sind unter 
$\Diff[\R^*]{X_0}^*$ irreduzibel, und die Algebra 
 \bdm
 \Diff[\R^*]{X_0}^*\ =\ \Aff{\C}{R,\, S,\, u\del_u,\, \del_{\theta}}
 \edm
ist nicht abelsch;
\item Die Vielfachheitenr"aume der $\SO(2,\R)$-Wirkung sind unter 
$\Diff[\SO(2,\R)]{X_0}^*$ zwar nicht irreduzibel, aber auch nicht zerlegbar,
keine zwei Vielfachheitenr"aume enthalten isomorphe Teildarstellungen
und die Algebra 
 \bdm
 \Diff[\SO(2,\R)]{X_0}^*\ =\ \Aff{\C}{u,\, u\del_u,\, \del_{\theta}}
 \edm
$\Diff[\SO(2,\R)]{X_0}^*$ ist nicht abelsch;
\item Die Vielfachheitenr"aume der $\SO(2,\R)$-Wirkung sind unter 
$\Diff[\SO(2,\R)]{X_0}$ irreduzibel, und die Algebra 
$\Diff[\SO(2,\R)]{X_0}$ ist nicht abelsch.
\end{enumerate}
\end{prop}
\begin{proof}
Mit Abbildung \ref{kegel-koordring} bleibt hier fast nichts zu 
beweisen, da sie die Operation von $\Diff{X_0}^*$ auf $A$ vollst"andig
beschreibt. Betrachten wir zun"achst die $\R^*$-Wirkung: weil es
keine nicht-konstanten invarianten Funktionen gibt, sind die
Vielfachheiten endlich, und die Irreduzibilit"at  der "`Spalten"' wird durch 
die invarianten Derivationen $R$ und $S$ bewerkstelligt. Invariante
Differentialoperatoren h"oherer Ordnung als die von den genannten erzeugten
kann es in $\Diff{X_0}^*$ nicht geben, weil es keine zur Multiplikation
mit den Koordinatenfunktionen "`inversen"' Derivationen gibt (etwa der 
Form $\del_u$); somit scheidet die M"oglichkeit, in Abbildung 
\ref{kegel-koordring} etwa erst nach rechts, dann mit $R$ oder $S$ herauf oder
herunter und dann in die urspr"ungliche Spalte zur"uck zugehen, aus.  
Wegen Gleichung \ref{rel-RS} h"atte man $\del_{\theta}$ in der Liste der 
erzeugenden Elemente auch weglassen k"onnen. 

F"ur die  $\SO(2,\R)$-Wirkung
ist $u$ eine invariante Funktion, deswegen sind die Vielfachheiten 
unendlich. Die einzigen $G$-invarianten
Derivationen sind Produkte aus $u$ mit $u\del_u$ und $\del_{\theta}$.
Alle aus $R$, $S$ und den Koordinatenfunktionen konstruierbaren "`Pfade"' in 
Abbildung \ref{kegel-koordring}, die in der gleichen Zeile anfangen und 
aufh"oren, lassen sich als Kombination einer Multiplikation mit einer 
$u$-Potenz und den Elementen $u\del_u$, $\del_{\theta}$ ausdr"ucken, deswegen 
liefern sie keine neuen $G$-invarianten  Differentialoperatoren.
Jeder Teilraum, der aus Elementen mit $u$-Grad gr"o"ser gleich einer
vorgegebenen nat"urlichen Zahl besteht, ist invariant; aus dem gleichen Grunde
sind die Vielfachheitenr"aume nicht zerlegbar. Zwei solche k"onnen keine
isomorphen Teildarstellungen haben, weil sich zwei Zeilen immer in ihrem
$\del_{\theta}$-Eigenwert unterscheiden.

Wirklich zu beweisen gibt es nun etwas bei der letzten Aussage. Die
Algebra $\Diff{X_0}$ aller algebraischer Differentialoperatoren des Kegels 
kann man nicht bestimmen; aber wir k"onnen einen $\SO(2,\R)$-invarianten
Operator hinschreiben, mit dem die Vielfachheitenr"aume irreduzibel werden.
Der Operator
 \bdm
 D\ :=\ \frac{1}{u}(u\del_u+i\del_{\theta})(u\del_u-i\del_{\theta})
 \edm
ist genau so gew"ahlt, da"s er den $u$-Grad um Eins absenkt, sofern er 
nicht auf ein Element der Gestalt $u^ne^{\pm i n \theta}$ wirkt, welches
er annuliert. Mit dem Kriterium aus Lemma \ref{diff-ops-krit} ist $D$ ein
zul"assiger Differentialoperator zweiter Ordnung und offensichtlich
invariant. In der Tat, in die Koordinaten $x$, $y$ und $z$ zur"uck"ubersetzt
schreibt sich $D$ als
 \bdm
 D\ =\ z(\del_x^2+\del_y^2)
 \edm
und verdient damit den Namen "`algebraisch"'.
\end{proof}
\end{exa}
\begin{NB}
Ein kurzer Blick auf den (glatten !) unendlich langen Zylinder $Z$ mit
Koordinatenring $\Aff{\C}{u, e^{i\theta}, e^{-i\theta} }$ erkl"art, warum
hier f"ur die  $S^1$-Wirkung Irreduzibilit"at vorliegt: in dieser
Situation ist n"amlich $\del_u$ eine $G$-invariante Derivation.
\end{NB}
\noindent
Dieses Beispiel ist noch relativ gutm"utig, und man k"onnte vermuten,
da"s Satz~\ref{isotyp-zerl-aff} vielleicht in einer schw"acheren 
Version (etwa Unzerlegbarkeit statt Irreduzibilit"at) gilt. Das
folgende Beispiel einer kubischen Regefl"ache wird zeigen, da"s sich diese 
Hoffnung  nicht erf"ullt.
\begin{figure}
\noindent
\centering
\begin{minipage}[b]{.46\linewidth}
\begin{turn}{-90}
\epsfig{figure=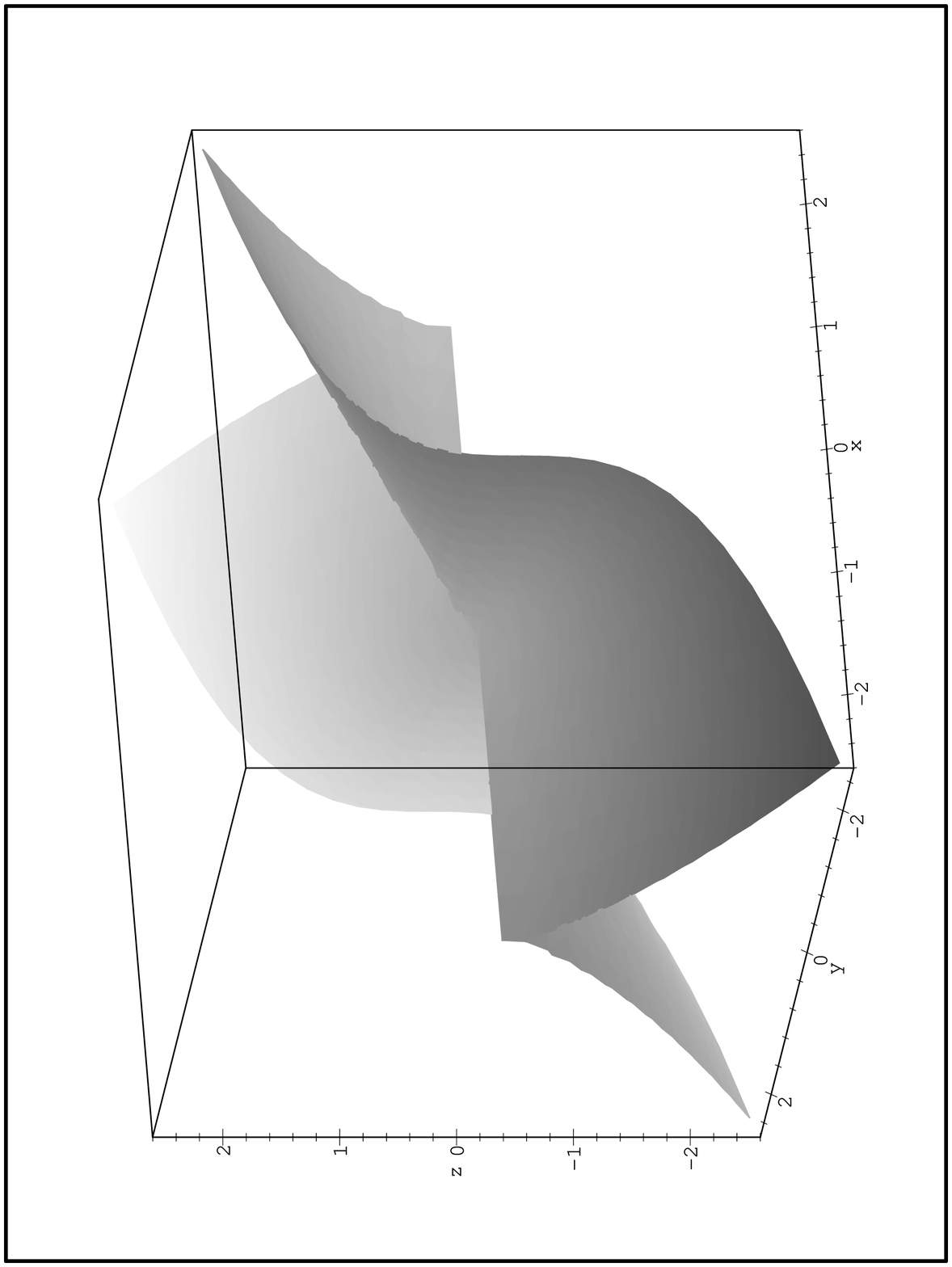, width=\linewidth, height=\linewidth}
\end{turn}
\caption{$z^3=xy^2$}\label{regel-1}
\end{minipage}
\end{figure}
%
\begin{exa}[Kubische Regelfl"ache $z^3=xy^2$]
Wir betrachten die irreduzible kubische Regelfl"ache im $\R^3$
(vgl.\ Abb.\ \ref{regel-1})
 \bdm
 Y_0\ :=\ \{(x,y,z)\in\R^3\ |\ z^3=xy^2\}\, . 
 \edm
Sie ist  eine reelle Form der komplexen Fl"ache
 \bdm
 Y \ :=\ \{(x,y,z)\in\C^3\ |\ z^3=xy^2\}\, , 
 \edm
und der singul"are Ort besteht in beiden F"allen aus der $x$-Achse.
Wir idenzitifieren deswegen wieder $\Aff{\C}{Y}$ mit $\Aff{\R}{Y_0}_{\C}$
und w"ahlen f"ur $Y_0$ die Parametrisierung
 \bdm
 x\ =\ rs^2 ,\quad y\ =\ r/s ,\quad z\ =\ r,
 \edm
erhalten also als regul"are Funktionen
 \bdm
 B\ :=\ \Aff{\R}{Y_0}_{\C}\ \cong\ \Aff{\C}{r,r/s,rs^2}\, ,
 \edm
die wir ebenfalls graphisch  dargestellt haben (Abb.\ 
\ref{regelflaeche-1-koordring}). Um die Wirkung der Algebra 
$\Diff{Y_0}^*$ hierauf zu beschreiben,
bestimmen wir wieder die Derivationen $\Der Y_0$. 
Den Koordinaten entsprechen zun"achst die Derivationen
 \bdm
 r\del_r\ =\ x\del_x+y\del_y+z\del_z\ \ \text{ und }\ \ 
 s\del_{s}\ =\ 2x\del_x-y\del_y\, .
 \edm
Sodann bestimmen wir mit der gleichen Methode wie in Lemma~\ref{deriv-kegel}
alle Derivationen von $B$, weswegen wir auf den Beweis verzichten:
\begin{lem}[Derivationen der Regelfl"ache $z^3=xy^2$]
\begin{enumerate}
\item Ein Element der Gestalt $q\del_r$ liegt genau dann in $\Der Y_0$,
wenn 
 \bdm
 q\in C\ :=\ B-\{\C\cdot r^n s^m\ |\ 2n=m\text { oder } n=-m\}\, ;
 \edm
\item Ein Element der Gestalt $q\del_s$ liegt genau dann in $\Der Y_0$,
wenn $\,q\in C$ oder $q\in s\cdot B$;
\item Ein Element der Gestalt $q_1\del_r + q_2\del_s$ liegt genau dann 
in $\Der Y_0$, wenn es Summe zweier Derivationen der gerade genannten
Gestalt oder in der $B$-linearen H"ulle der Derivation
 \bdm
 D\ :=\ \frac{r^2}{s^2}\del_r+\frac{r}{s}\del_s
 \edm
enthalten ist.
\end{enumerate}\qed
\end{lem}
\begin{figure}
\begin{diagram}[height=2.8em,width=1.4em,labelstyle=\scriptstyle,postscript,tight]
  & & & & & & & & & & & & r^3s^6 & & \rTo& & r^4s^6&\rDots & & & 
    \hphantom{x} & \rDots & & & \hphantom{x} \\
  & & & & & & & & & & & \ruTo(4,4) & & \rdLine(2,2)^D & & & & \rdLine(2,2)^D 
     & & & & & & \\
  & & & & & & & & & & & & 0 & & \rdTo(3,3) & & r^4s^5 & \rDots & & & 
    \hphantom{x} & \rDots & & & \hphantom{x} \\
  & & & & & & & & & & & & & & & & & \rdLine(2,2)^D & & \rdTo(2,2) & & & & & \\
  & & & & & & & & r^2 s^4 & & \rTo & & r^3 s^4 & & \rTo& & r^4s^4 &\rDots & & &
     \hphantom{x} & & & & \\
  & & & & & & &  \ruTo(4,4)& & \rdLine(2,2)^D & & & & \rdLine(2,2)^D & & & & 
     \rdLine(2,2)^D & & \rdTo(2,2) & & & & \\
  & & & & & & & & 0 & & \rdTo(3,3) & & r^3 s^3& & \rTo& & r^4s^3 &\rDots & & & 
     \hphantom{x} & & & & \\
  & & & & & & & & & & & & & \rdLine(2,2)^D & & \rdTo(2,2) & & & & \rdTo(2,2) &
     & & & \\
  & & & & rs^2& & \rTo & & r^2 s^2& & \rTo & & r^3 s^2 & & \rTo& & r^4s^2 &
     \rDots & & & \hphantom{x} & & & & \\
  & & & \ruTo(4,4)& & \rdLine(2,2)^D & & & & \rdLine(2,2)^D & & & & 
     \rdLine(2,2)^D & & \rdTo(2,2) & & & & & & & & \\
  & & & & 0 & &\rdTo(3,3)& & r^2 s & & \rTo & & r^3 s & & \rTo& & r^4 s & &
     & & & & & & \\
  & & & & & & & & & \rdLine(2,2)^D & & \rdTo(2,2) & & \rdLine(2,2)^D & & 
     \rdTo(2,2) & & & & & & & & \\
1 & & \rTo& & r & & \rTo & & r^2 & & \rTo & & r^3 & & \rTo& & r^4 & & & & &
     & & & \\
  & \rdTo(4,2)& & & & \rdLine(2,2)^D & & & & \rdLine(2,2)^D & & \rdTo(2,2) & & 
    \rdLine(2,2)^D & & \rdTo(2,2)& & & & & & & & \\
  & & & & r/s & & \rTo & & r^2/ s & & \rTo & & r^3/s & & \rTo& & r^4/s & &
     & & & & & \\
  & & & & & \rdTo(4,2)\rdLine(3,3) & & \rdTo(2,2)& & \rdLine(2,2)^D & & 
    \rdTo(2,2) & & \rdLine(2,2)^D & &\rdTo(2,2) & & & & & & & & \\
  & & & & 0 & & & & r^2/s^2& & \rTo & & r^3/s^2& & \rTo& & r^4/s^2 & \rDots &
     & &  \hphantom{x} & & & & \\
  & & & & & & & \rdTo(2,2)_D& & \rdTo(4,2)\rdLine(3,3)& & \rdTo(2,2)& & 
      \rdLine(2,2)^D & & \rdTo(2,2)& & \rdLine(2,2)^D & & & & & & \\
  & & & & & & & & 0 & & & & r^3/s^3& & \rTo& & r^4/s^3 & \rDots & & & 
      \hphantom{x} & & & & \\
  & & & & & & & & & & & \rdTo(2,2)_D& &  \rdTo(4,2)\rdLine(3,3)& &\rdTo(2,2) &
     & \rdLine(2,2)^D & &\rdTo(2,2) & & & & \\
  & & & & & & & & & & & & 0 & & & & r^4/s^4 &\rTo & & & \hphantom{x} &
      \rDots & & & \hphantom{x} \\
  & & & & & & & & & & & & & & & \rdTo(2,2)_D& & \rdTo(4,2) & & \rdTo(2,2)& &
    & &\\
  & & & & & & & & & & & & &  & & & 0 & & & & \hphantom{x}& \rDots& & &
    \hphantom{x}\\
\end{diagram}\caption{Koordinatenring der Regelfl"ache  $z^3=xy^2$}
\label{regelflaeche-1-koordring}
\end{figure}
%
\noindent
Insbesondere gibt es weder den $r$- noch den $s$-Grad absenkende Derivationen,
und nur die relativ uninteressanten Elemente  $r\del_r$ und $s\del_s$ erhalten
den Grad.
Die Wirkung von $D$ ist (bis auf von Null verschiedene Vielfache)
in Abbildung \ref{regelflaeche-1-koordring} dargestellt; im Gegensatz zu
zu den Elementen $R$ und $S$, die beim Kegel aufgetreten waren, ver"andert
$D$ auch den "'Spaltengrad"'. 
Wegen der geringeren Symmetrie haben wir
als Gruppenwirkung hier nur die Wahl zwischen der multiplikativen
$\R^*$-Wirkung oder der der $s$-Koordinate entsprechenden $\R^*$-Wirkung
in jeder $x$-$y$-Ebene.
%
\begin{prop}
\begin{enumerate}
\item Die Vielfachheitenr"aume der multiplikativen $\R^*$-Wirkung zerfallen
unter der Wirkung der abelschen Algebra 
 \bdm
 \Diff[\R^*]{Y_0}^*\ =\ \Aff{\C}{r\del_r,\, s\del_s}
 \edm
in eindimensionale irreduzible Darstellungen;
\item Die Vielfachheitenr"aume der $s$-Wirkung sind unter 
$\Diff[s]{Y_0}^*$ weder irreduzibel noch zerlegbar, und
keine zwei Vielfachheitenr"aume enthalten isomorphe Teildarstellungen
der Algebra 
 \bdm
 \Diff[s]{Y_0}^*\ =\ \Aff{\C}{r,\, r\del_r,\, s\del_s}\, .
 \edm
\end{enumerate}
\end{prop}
\begin{proof}
Die beiden Behauptungen "uberlegt man sich leicht anhand von
Abbildung  \ref{regelflaeche-1-koordring} und dem vorangegangenen Lemma.
Wir bemerken kurz, da"s man hier nach dem gleichen Schema wie f"ur den
Kegel invariante Differentialoperatoren hinschreiben kann, die nicht
in $\Diff{Y_0}^*$ liegen; zum Beispiel  ist
 \bdm
 T\ :=\ \frac{1}{s}(2r\del_r-s\del_s)(r\del_r +s\del_s)\ =\
 9z^2\del^2_{xy}+3y^2\del^2_{yz}+6xy\del^2_{xz}+2y\del_y+2zy\del^2_{zz}
 \edm
ein Operator, welcher unter der skalaren $\R^*$-Wirkung invariant ist, den
$s$-Grad um Eins absenkt und die "`L"ucken"' im Koordinatenring gerade auf Null
abbildet.
\end{proof}
Die erste Aussage ist ein Gegenbeispiel zu Proposition \ref{multfree-diffops}
im singul"aren Fall: die Algebra $\Diff[G]{M}^*$ ist hier abelsch, obwohl
die Wirkung \emph{nicht} vielfachheitenfrei ist.
\end{exa}
%
\begin{thank}
Mein herzlichster Dank gilt jenen, die mich bei der
Anfertigung dieser Arbeit mit wertvollen Hinweisen, Diskussionen
oder auch nur ihrer Geduld begleitet haben: dies sind
Thomas Friedrich (Humboldt-Universit"at zu Berlin), Roe Goodman
(Rutgers University) und Karl-Hermann Neeb (TU Darmstadt). 
\end{thank}

%

\end{document}